\documentclass[11pt,a4paper]{article}
\date{}
\makeatletter
\@addtoreset{equation}{section}
\makeatother

\usepackage[colorlinks,
            linkcolor=blue,
            anchorcolor=blue,
            citecolor=blue]{hyperref}
\usepackage{hyperref}
\hypersetup{colorlinks=true,
            linkcolor=blue,
            anchorcolor=blue,
            citecolor=blue}
\usepackage{hyperref,amsmath,amssymb,amscd}
\usepackage{amssymb,amsmath,enumerate}
\usepackage{indentfirst}
\usepackage{url}
\usepackage[capitalise]{cleveref}
\usepackage{amsthm}
\usepackage{amsmath}
\usepackage{diagbox}

\usepackage{booktabs}
\usepackage{diagbox}
\usepackage{makecell}

\AtBeginDocument{\hypersetup{pdfborder={0 0 0.01}}}
\newtheorem{theorema}{Theorem}[section]

\newtheorem{thm}{Theorem}[section]

\newtheorem{lem}[thm]{Lemma}
\newtheorem{prop}[thm]{Proposition}
\theoremstyle{definition}

\newtheorem{rem}[thm]{\bf Remark}
\allowdisplaybreaks[4]

\usepackage{enumitem}
\usepackage{extarrows}
\usepackage{amsthm}
\usepackage{setspace}
\hypersetup{urlcolor=blue, citecolor=red}
\allowdisplaybreaks[4]
\usepackage{cite}

\begin{document}

\title{\bf
On the extremal functions of second order uncertainty principles: symmetry and symmetry breaking\footnote{This work was supported by National Natural Science Foundation of China (No. 12371120).}}
\author{{Xiao-Ping Chen,\ \ Chun-Lei Tang\footnote{Corresponding author.\newline
\indent\,\,\, \emph{E-mail address:}   xpchen\_maths@163.com (X.-P. Chen); tangcl@swu.edu.cn (C.-L. Tang).}}\\
{\small \emph{School of Mathematics and Statistics, Southwest University,  Chongqing {\rm400715},}}\\
{\small \emph{People's Republic of China}}}
\maketitle
\baselineskip 17pt

\noindent {\bf Abstract}:
\ This paper focus on the symmetry and symmetry breaking about the second order Hydrogen Uncertainty Principle. \emph{Firstly}, by choosing a suitable test function, we give a negative answer to the conjecture presented by Cazacu, Flynn and Lam in [\emph{J. Funct. Anal.} \textbf{283} (2022), Paper No. 109659, 37 pp] for $N\in\{2,3\}$, and emphasizing the symmetry breaking phenomenon. \emph{Secondly}, we obtain a family of sharp weighted second order Hydrogen Uncertainty Principle, and prove the extremal functions are radial, which extends the work of Duong and Nguyen [The sharp second order Caffareli-Kohn-Nirenberg inequality and stability estimates for the sharp second order uncertainty principle,  arXiv:2102.01425].

\vspace{0.25em}


\vspace{0.25em}


\section{Introduction and main results}

\noindent This paper aims to study second order Hydrogen Uncertainty Principles (HyUP for short), which can be stated more precisely as follows.
\begin{enumerate}
[itemsep=0pt,
topsep=2pt,
parsep=2pt]

\item[($a$)]
Firstly, we consider the following second order HyUP,
\begin{equation*}
\int_{\mathbb{R}^N}|\Delta u|^2\mathrm{d}x
\int_{\mathbb{R}^N}\left|\nabla u\right|^2\mathrm{d}x
\ge C(N)
\left(\int_{\mathbb{R}^N}
\frac{\left|\nabla u\right|^2}{|x|}
\mathrm{d}x\right)^2,\quad \forall u\in C_c^\infty(\mathbb{R}^N),
\end{equation*}
which has been proved by Cazacu, Flynn and Lam in \cite[Theorem 2.3]{Cazacu22} for higher dimensions $N\geq 5$ with $C(N)=\frac{(N+1)^2}{4}$. Moreover, the authors in \cite[Remark 2.4]{Cazacu22} conjectured that
$C(N)=\frac{(N+1)^2}{4}$ for $2\leq N\leq 4$ . In Theorem \ref{thm-3} below, by choosing a suitable test function, we will give a negative answer to this conjecture, i.e., $C(N)<\frac{(N+1)^2}{4}$ for $N\in\{2,3\}$, and emphasize the symmetry breaking phenomenon.

\item[($b$)]
Secondly, by using spherical harmonics decomposition, we study the following sharp weighted second order HyUP:
\begin{equation}\label{8.16-4}
\int_{\mathbb{R}^N}\frac{|\Delta u|^2}{|x|^{2\alpha}}\mathrm{d}x
\int_{\mathbb{R}^N}\left|\nabla u\right|^2\mathrm{d}x
\ge C(N,\alpha)
\left(\int_{\mathbb{R}^N}
\frac{\left|\nabla u\right|^2}{|x|^{\alpha+1}}
\mathrm{d}x\right)^2,\quad \forall u\in
C_c^\infty(\mathbb{R}^N),
\end{equation}
where
\begin{equation*}
C(N,\alpha)=
\begin{cases}
\frac{(N+3\alpha+1)^2}{4},
\ \ &
\mbox{if}
\
N\ge2,
\
\alpha>-1
\
\mbox{and}
\
N\ge5\alpha+5,\\[1.5mm]
\frac{(3\alpha+5)^2}{4},
\ \ &
\mbox{if}
\
N=4
\
\mbox{and}
\
-1<\alpha\le0;
\end{cases}
\end{equation*}
we also obtain the sufficient and necessary condition such that the equality in \eqref{8.16-4} holds, see Theorems \ref{thmwsupsc-2} and \ref{thmwsupscN4} below. These results improve and generalize that of \cite{Cazacu22}.

\end{enumerate}

The rest of this section is organized as follows. In Section \ref{sect-1.1}, we recall some significant results of the first order and second order HyUP (or some related inequalities). Our main results and related comments are described in Section \ref{sect-1.2}.

\subsection{Overview and motivation}\label{sect-1.1}

\noindent We begin with the first order HyUP that has important applications in mathematics and  physics (for example, which can be used to study the stability of a hydrogenic atom in a magnetic field, see e.g., \cite{Frohlich86}): for $N\ge2$,
\begin{equation}\label{8.11-8}
\int_{\mathbb{R}^N}
|\nabla u|^2\mathrm{d}x
\int_{\mathbb{R}^N}
\left|u\right|^2\mathrm{d}x
\ge\frac{(N-1)^2}{4}
\left(\int_{\mathbb{R}^N}
\frac{\left|u\right|^2}{|x|}
\mathrm{d}x\right)^2,
\end{equation}
where the constant $\frac{(N-1)^2}{4}$ in \eqref{8.11-8} is sharp and can be attained by $u(x)=ae^{-b|x|}$ with $a\in\mathbb{R}$ and $b>0$, see \cite[Example 1.6]{Frank11} for details.

The HyUP \eqref{8.11-8} is a special case of the first order interpolation inequalities which is known as   Caffarelli-Kohn-Nirenberg (CKN for short) inequalities (see \cite{Caffarelli84}). A large number of scholars have paid attention to CKN inequalities in the past several decades due to their applications in many area of mathematics and physics, we refer the readers to see \cite{Catrina09,Catrina01,Cazacu21,
Cazacu24-2,Costa08,Dong18,
Xia07,Do23,Duy22,Frank08,
Lam20,Vazquez00,DelPino03,
Nguyen15,Talenti76,Brezis18,
Cordero-Erausquin04,DelPino02} for first order CKN inequalities and the references therein.

Compared with the first order CKN inequalities, there exist few results about the higher order CKN inequalities except \cite{Cazacu22,Cazacu23,
Duong25-2,Duong25,Lin86}. We first introduce some basic notations in order to state the previous results about higher order CKN inequalities. For $\alpha,\beta$ satisfying $N-2\alpha>0$ and $N-2\beta>0$, let $\mathcal{H}^2_{\alpha,\beta}
(\mathbb{R}^N)$ be the second order Sobolev space which is the closure of $C_c^\infty(\mathbb{R}^N)$ with respect to the norm
\begin{equation*}
\|u\|_{\mathcal{H}^2_{\alpha,\beta}
(\mathbb{R}^N)}
:=\left(\int_{\mathbb{R}^N}
\frac{|\Delta u|^2}
{|x|^{2\alpha}}\mathrm{d}x
+\int_{\mathbb{R}^N}
\frac{|\nabla u|^2}
{|x|^{2\beta}}\mathrm{d}x
\right)^{\frac{1}{2}}.
\end{equation*}

In the recent paper \cite{Cazacu22}, Cazacu, Flynn and Lam considered the following  second order HyUP in higher dimensions $N\ge5$, and proved the following results.

\begin{theorema}
[{\!\rm{\!\cite[Theorem 2.3]{Cazacu22}}}]
\label{thm-c}
Let $N\ge5$. Then the following inequality holds
\begin{equation}\label{sUP}
\int_{\mathbb{R}^N}
|\Delta u|^2
\mathrm{d}x
\int_{\mathbb{R}^N}
|\nabla u|^2\mathrm{d}x
\geq C(N)
\left(\int_{\mathbb{R}^N}
\frac{|\nabla u|^2}{|x|}
\mathrm{d}x\right)^2,
\end{equation}
for any $u\in
\mathcal{H}^2_{0,0}(\mathbb{R}^N)$, where $C(N)=\frac{(N+1)^2}{4}$ is the sharp constant of \eqref{sUP} and can be attained by $u(x)=a(1+b|x|)\exp(-b|x|)$ with   $a\in\mathbb{R}$ and $b>0$.
\end{theorema}

Note that their main idea is to apply spherical harmonics decomposition, but their proof fails for lower dimensions $2\le N\le4$. Moreover, let us denote $\mathcal{R}_1:
=\frac{x}{|x|}\cdot\nabla$ and $\mathcal{R}_2:
=\mathcal{R}_1^2
+\frac{N-1}{r}\mathcal{R}_1$, where $\mathcal{R}_2:
=\mathcal{R}_1\circ\mathcal{R}_1$, which represent the radial derivative and the radial Laplacian, respectively. Then, in \cite[Theorem 2.4]{Cazacu22}, the authors also proved that, for all  $N\geq 2$,
\begin{equation*}
\int_{\mathbb{R}^N}
|\mathcal{R}_2 u|^2\mathrm{d}x
\int_{\mathbb{R}^N}
\left|\mathcal{R}_1 u\right|^2 \mathrm{d}x
\ge \frac{(N+1)^2}{4}
\left(\int_{\mathbb{R}^N}
\frac{\left|\mathcal{R}_1 u\right|^2}{|x|}
\mathrm{d}x\right)^2,
\end{equation*}
and the equality holds by $u(x)=a(1+b|x|)e^{-b|x|}$ with  $a\in\mathbb{R}$ and $b>0$. That is to say, for any radial function $u\in
\mathcal{H}^2_{0,0}(\mathbb{R}^N)$, $C(N)=\frac{(N+1)^2}{4}$ for $2\le N\le4$. Therefore, for all $u\in
\mathcal{H}^2_{0,0}(\mathbb{R}^N)$ (without radiality assumption), $C(N)\le\frac{(N+1)^2}{4}$ for $2\le N\le4$. \emph{In} \cite[Remark 2.4]{Cazacu22}, \emph{the authors conjectured that the sharp constant in \eqref{sUP} satisfies $C(N)
=\frac{(N+1)^2}{4}$ for $2\leq N\leq 4$}. We will give a negative answer to this conjecture for $N\in\{2,3\}$ in Theorem \ref{thm-3} below.

It is worthy to mention that the authors in \cite{Cazacu22} also studied the so-called Heisenberg Uncertainty Principle, and obtained the following results.

\begin{theorema}
[{\!\rm{\!\cite[Theorem 2.1]{Cazacu22}}}]
\label{thm-d}
Let $N\ge1$. Then for any $u\in
\mathcal{H}^2_{0,-1}(\mathbb{R}^N)$,
\begin{equation}\label{sUP-2}
\int_{\mathbb{R}^N}
|\Delta u|^2
\mathrm{d}x
\int_{\mathbb{R}^N}
|x|^2|\nabla u|^2\mathrm{d}x
\geq \frac{(N+2)^2}{4}
\left(\int_{\mathbb{R}^N}
|\nabla u|^2\mathrm{d}x\right)^2,
\end{equation}
where $\frac{(N+2)^2}{4}$ is the sharp constant and is attained by $u(x)=a\exp(-b|x|^2)$,  $a\in\mathbb{R}$ and $b>0$.
\end{theorema}

Subsequently, Cazacu, Flynn and Lam in \cite{Cazacu23} improved \eqref{sUP} and \eqref{sUP-2} to the following  weighted version.

\begin{theorema}
[{\!\rm{\!\cite[Corollary 2.4]{Cazacu23}}}]
\label{thm-b}
Let $N\ge5$. The following inequality
\begin{equation}\label{8.11-6}
\int_{\mathbb{R}^N}
|\Delta u|^2
\mathrm{d}x
\int_{\mathbb{R}^N}
\frac{|\nabla u|^2}
{|x|^{2\beta}}\mathrm{d}x
\geq C(N,0,\beta)
\left(\int_{\mathbb{R}^N}
\frac{|\nabla u|^2}
{|x|^{\beta+1}}
\mathrm{d}x\right)^2
\end{equation}
holds for any $u\in \mathcal{H}^2_{0,\beta}(\mathbb{R}^N)$. Furthermore,
\begin{enumerate}
[itemsep=0pt,
topsep=2pt,
parsep=2pt]

\item[(1)]
if $\beta<1$, then $C(N,0,\beta)
=\frac{(N-\beta+1)^2}{4}$ is the sharp constant of \eqref{8.11-6} and can be attained by the functions $u$ satisfying $\nabla u=a\exp\left(-\frac{b}{1-\beta}
|x|^{1-\beta}\right)x$ with  $a\in\mathbb{R}$ and $b>0$;

\item[(2)]
if $\beta>1$, then $C(N,0,\beta)
=\frac{(N+\beta-1)^2}{4}$ is the sharp constant of \eqref{8.11-6} and can be attained by the functions $u$ satisfying $\nabla u=a|x|^{-N} \exp\left(-\frac{b}{1-\beta}
|x|^{1-\beta}\right)x$ with $a\in\mathbb{R}$ and $b<0$.

\end{enumerate}
\end{theorema}

In particular, when $\beta=1$ in \eqref{8.11-6}, it becomes the well-known sharp Hardy-Rellich inequality, see \cite{Tertikas07,Beckner08,Cazacu20} for details.

Along the same line of thought, the authors in \cite{Cazacu23} also obtained the following sharp weighted second order Heisenberg Uncertainty Principle.

\begin{theorema}
[{\!\rm{\!\cite[Corollary 2.5]{Cazacu23}}}]
\label{thm-a}
Let $N\ge2$ and $\alpha\in\mathbb{R}$. There holds
\begin{equation}\label{gwsup-1}
\int_{\mathbb{R}^N}
\!\frac{|\Delta u|^2}
{|x|^{2\alpha}}
\mathrm{d}x
\int_{\mathbb{R}^N}
|x|^{2\alpha+2}|\nabla u|^2\mathrm{d}x
\geq \frac{(N+4\alpha+2)^2}{4}
\left(\int_{\mathbb{R}^N}
|\nabla u|^2
\mathrm{d}x\right)^2,
\end{equation}
for any $u\in
\mathcal{H}^2_{\alpha,-\alpha-1}
(\mathbb{R}^N)$. Moreover, if either $\alpha+1>0$ or $N+4\alpha+2>0$, then $\frac{(N+4\alpha+2)^2}{4}$ is the sharp constant of \eqref{gwsup-1} and can be attained by $u(x)=a
\exp\left[-b|x|^{2(\alpha+1)}\right]$, $a\in\mathbb{R}$ and $b>0$.
\end{theorema}

Recently, Duong and Nguyen in \cite{Duong25} established a family of more general second order CKN inequalities: let $N\geq 1$ and $\alpha,\beta$ be such that
\begin{equation*}
N-2\alpha>0,
\
N-2\beta>0,
\
N-\alpha-\beta-1>0,
\end{equation*}
then
\begin{equation}\label{gwsup}
\int_{\mathbb{R}^N}
\frac{|\Delta u|^2}
{|x|^{2\alpha}}
\mathrm{d}x
\int_{\mathbb{R}^N}
\frac{|\nabla u|^2}
{|x|^{2\beta}}\mathrm{d}x
\geq C(N,\alpha,\beta)
\left(\int_{\mathbb{R}^N}
\frac{|\nabla u|^2}
{|x|^{\alpha+\beta+1}}
\mathrm{d}x\right)^2,
\end{equation}
for any $u\in \mathcal{H}^2_{\alpha,\beta}
(\mathbb{R}^N)$. Furthermore, if $\alpha-\beta+1>0$ and $N+2\alpha>0$, the constant $C(N,\alpha,\beta)$ in \eqref{gwsup} is sharp and  satisfies
\begin{equation*}
C(N,\alpha,\beta)
\geq\inf_{k\in\mathbb{N}^*\cup\{0\}}
\dfrac{1+\min\left\{0,
\frac{8\beta k}
{\left(N+2k-2\beta-2\right)^2}
\right\}}
{\left[1+\max\left\{0,
\frac{4(\alpha+\beta+1)k}
{\left(N+2k-\alpha-\beta-3\right)^2
}\right\}\right]^2}
\left(\frac{N+2k+3\alpha-\beta+1}
{2}\right)^2.
\end{equation*}
However, the authors did not give the explicit form of $C(N,\alpha,\beta)$ and its extremal functions.

Note that \eqref{gwsup} contains \eqref{sUP}, \eqref{sUP-2}, \eqref{8.11-6} and \eqref{gwsup-1}. For the sake of comparison, we present the previous results of \eqref{gwsup} about sharp constants and extremal functions in the table below:

\vspace*{0.25em}

\begin{center}
\scriptsize
\begin{tabular}
{l|c|c|c|c}
  \toprule
  \diagbox
  {$\alpha,\beta$}
  {results}
  {$N$}
  &
  $N=2$
  &
  $N=3$
  &
  $N=4$
  &
  $N\ge5$
  \\
  \hline
  $\alpha=0$, $\beta=0$
  &

  &

  &

  &
  \cite[Theorem 2.3]{Cazacu22}
  \\
  \hline
  $\alpha=0$, $\beta=1$
  &
  $\diagdown$
  &
  \makecell[c]
  {\cite[Corollary 1]{Beckner08}
  \\\cite[Theorem 1.1]{Cazacu20}}
  &
  \makecell[c]
  {\cite[Corollary 1]{Beckner08}
  \\\cite[Theorem 1.1]{Cazacu20}}
  &
  \makecell[c]
  {\cite[(1.8)]{Tertikas07}
  \\\cite[Corollary 1]{Beckner08}
  \\\cite[Theorem 1.1]{Cazacu20}}
  \\
  \hline
  $\alpha=0$, $\beta=-1$
  &
  \cite[Theorem 2.1]{Cazacu22}
  &
  \cite[Theorem 2.1]{Cazacu22}
  &
  \cite[Theorem 2.1]{Cazacu22}
  &
  \cite[Theorem 2.1]{Cazacu22}
  \\
  \hline
  $\alpha=0$, $\beta\in\mathbb{R}$
  &

  &

  &

  &
  \cite[Corollary 2.4]{Cazacu23}
  \\
  \hline
  $\alpha\in\mathbb{R}$, $\beta=-\alpha-1$
  &
  \cite[Corollary 2.5]{Cazacu23}
  &
  \cite[Corollary 2.5]{Cazacu23}
  &
  \cite[Corollary 2.5]{Cazacu23}
  &
  \cite[Corollary 2.5]{Cazacu23}
  \\
  \hline
  $\alpha\in\mathbb{R}$, $\beta=\alpha+1$
  &
  \cite[Corollary 1.3]{Cassano22}
  &
  \cite[Corollary 1.3]{Cassano22}
  &
  \cite[Corollary 1.3]{Cassano22}
  &
  \cite[Corollary 1.3]{Cassano22}
  \\
  \bottomrule
\end{tabular}
\end{center}

\vspace*{0.5em}

Inspired by the papers mentioned above, it is natural to ask a question that \emph{when $\alpha\in\mathbb{R}$ and $\beta=0$ in \eqref{gwsup}, what is the explicit form of $C(N,\alpha,0)$ and its extremal functions}. We will give a positive answer to it in Theorems \ref{thmwsupsc-2} and \ref{thmwsupscN4} below.

\subsection{Main results and some comments}\label{sect-1.2}

\noindent The main idea in this paper is to adapt the technique of decomposing a function into spherical harmonicas, which is a useful method (see \cite{Cazacu20,Cazacu22,
Tertikas07,Vazquez00} and the references therein). For $N\ge2$, we apply the coordinate transformation $x\in\mathbb{R}^N
\mapsto(r,\sigma)
\in(0,\infty)\times
\mathbb{S}^{N-1}$, and then for each $u\in C^\infty_c(\mathbb{R}^N)$,
\begin{equation}\label{1020-1}
u(x)=u(r\sigma)
=\sum^\infty_{k=0}
u_k(r)\phi_k(\sigma),
\end{equation}
where $\phi_k$ are the orthonormal eigenfunctions of the Laplace-Beltrami operator on $\mathbb{S}^{N-1}$ with the corresponding eigenvalue $c_k=k(N+k-2)$, namely, $-\Delta_{\mathbb{S}^{N-1}}\phi_k
=c_k\phi_k$ with $k\in\mathbb{N}^*\cup\{0\}$. Moreover, to simplify some calculations below, we will choose such that
\[
\int_{\mathbb{S}^{N-1}}
|\phi_k(\sigma)|^2\mathrm{d}\sigma
=\frac{|\mathbb{S}^{N+2k-1}|}
{|\mathbb{S}^{N-1}|},
\]
where $\mathrm{d}\sigma$ is the normalized surface measure on $\mathbb{S}^{N-1}$.

The Fourier coefficients $\{u_k\}_k$ belong to
$C_c^\infty([0,\infty))$ and satisfy $u_{k}(r)=O(r^k)$ as $r\to 0$, which allows us to consider the change of variables
\begin{equation*}
u_{k}(r)=r^kv_k(r),
\end{equation*}
where $v_k\in C^\infty_c([0,\infty))$ (see \cite[page 418]{Tertikas07} for more details). Therefore, the function $u$ can be represented by the following form
\begin{equation}\label{1031-2}
u(x)=u(r\sigma)
=\sum_{k=0}^{\infty}u_{k}(r)
\phi_{k}(\sigma)
=\sum_{k=0}^{\infty}r^{k}v_{k}(r)
\phi_{k}(\sigma).
\end{equation}
Then, direct calculations together with \cite[(2.3)--(2.5)]{Duong25-2}  imply the following identities
\begin{align*}
\int_{\mathbb{R}^{N}}
|\Delta u|^{2}\mathrm{d}x
&=\sum_{k=0}^{\infty}
\int_{\mathbb{R}^{N+2k}}
|\Delta v_k(|x|)|^{2}\mathrm{d}x,
\\
\int_{\mathbb{R}^{N}}
|\nabla u|^{2}\mathrm{d}x
&=\sum_{k=0}^{\infty}
\int_{\mathbb{R}^{N+2k}}
\left|\nabla v_{k}(|x|)\right|^{2}
\mathrm{d}x,\\
\int_{\mathbb{R}^{N}}
\frac{|\nabla u|^{2}}{|x|}
\mathrm{d}x
&=
\sum_{k=0}^{\infty}
\int_{\mathbb{R}^{N+2k}}
\left[
\frac{|\nabla v_k(|x|)|^{2}}{|x|}
+k
\frac{|v_k(|x|)|^{2}}{|x|^3}
\right]\mathrm{d}x.
\end{align*}

For $N\ge2$ and  $k\in\mathbb{N}^*\cup\{0\}$, let
\begin{align*}
C(N,k)
:=\inf_{u\in
W^{2,2}_{\mathrm{rad}}
(\mathbb{R}^{N+2k})
\setminus\{0\}}
\frac{\int_{\mathbb{R}^{N+2k}}
|\Delta u|^{2}\mathrm{d}x
\int_{\mathbb{R}^{N+2k}}
|\nabla u|^{2}\mathrm{d}x}
{\left[\int_{\mathbb{R}^{N+2k}}
\left(\frac{|\nabla u|^{2}}{|x|}
+k\frac{|u|^{2}}{|x|^3}\right)
\mathrm{d}x\right]^2}.
\end{align*}
From \eqref{sUP}, we see that  $C(N,0)=\frac{(N+1)^2}{4}$. It is easy to check that
\begin{equation*}
\inf_{k\in\mathbb{N}^*\cup\{0\}}
C(N,k)
=C(N).
\end{equation*}

\begin{thm}\label{thm-3}
For $N\in\{2,3\}$, there holds
\[
C(N,k)\ge
\frac{(N+2k+1)^2}{4\left[1+\frac{4k}
{(N+2k-3)^2}\right]^2},
\]
for all $k\in\mathbb{N}^*$. As a consequence, the sharp constant $C(N)$ of \eqref{sUP} satisfies
\begin{equation*}
\frac{(N-1)^4(N+3)^2}
{4[(N-1)^2+4]^2}
\le C(N,1)=C(N)
\leq\frac{N(N+4)(N^2-1)^2}
{4(N^2-N+4)^2}
<\frac{(N+1)^2}{4}.
\end{equation*}
\end{thm}

\begin{rem}
\rm
From \cite[Theorem 1.3]{Huang25}, we know that the extremal function of \eqref{sUP} exists, this together with $C(N)=C(N,1)$ implies that the extremal function of \eqref{sUP} must be of the following non-radial form
\[
u(x)=|x|v_1(|x|)\phi_1(\sigma),
\]
where $v_1(|x|)$ is radial and $\phi_1(\sigma)$ is defined by \eqref{1020-1}. This can be called symmetry breaking phenomenon.
\end{rem}

\begin{thm}\label{thmwsupsc-2}
Assume that $N\ge2$ and $\alpha>-1$ satisfy $N\geq 5\alpha+5$. For any $u\in \mathcal{H}^2_{\alpha,0}(\mathbb{R}^N)$,
there holds:
\begin{equation}\label{wsUP}
\int_{\mathbb{R}^N}
\frac{|\Delta u|^2}
{|x|^{2\alpha}}
\mathrm{d}x
\int_{\mathbb{R}^N}\left|\nabla u\right|^2\mathrm{d}x
\ge\frac{(N+3\alpha+1)^2}{4}
\left(\int_{\mathbb{R}^N}
\frac{\left|\nabla u\right|^2}
{|x|^{\alpha+1}}
\mathrm{d}x\right)^2,
\end{equation}
where the constant $\frac{(N+3\alpha+1)^2}{4}$
is sharp and is attained if and only if
\begin{equation*}
u(x)=a\left(1+b|x|^{\alpha+1}\right)
\newline
\exp\left(-b |x|^{\alpha+1}\right),
\ \
a\in\mathbb{R},
\
b>0.
\end{equation*}
\end{thm}

\begin{rem}
Compared with \cite[Theorem 2.3]{Cazacu22}
(or see Theorem \ref{thm-c} above), we
obtain a weighted version of the second
order HyUP \eqref{sUP} in
Theorem \ref{thmwsupsc-2};
more precisely, if $\alpha=0$ in
Theorem \ref{thmwsupsc-2}, it
reduces to \cite[Theorem 2.3]{Cazacu22} or the case $N\ge5$ of \cite[Theorem 1.2]{Huang25}.
\end{rem}

Notice that Theorem \ref{thmwsupsc-2} does not include the result of Huang and Ye \cite{Huang25} about the HyUP \eqref{sUP} in dimension $4$. Now, we are going to establish the weighted HyUP in $\mathbb{R}^4$.

\begin{thm}\label{thmwsupscN4}
Assume that $-1<\alpha\leq 0$.
For any $u\in \mathcal{H}^2_{\alpha,0}(\mathbb{R}^4)$, there holds
\begin{equation}\label{wsUPN4}
\int_{\mathbb{R}^4}
\frac{|\Delta u|^2}
{|x|^{2\alpha}}
\mathrm{d}x
\int_{\mathbb{R}^4}\left|\nabla u\right|^2\mathrm{d}x
\ge\frac{(3\alpha+5)^2}{4}
\left(\int_{\mathbb{R}^4}
\frac{\left|\nabla u\right|^2}
{|x|^{\alpha+1}}
\mathrm{d}x\right)^2,
\end{equation}
where the constant $\frac{(3\alpha+5)^2}{4}$
is sharp and is attained if and only if
\begin{equation*}
u(x)=a\left(1+b|x|^{\alpha+1}\right)
\exp\left(-b |x|^{\alpha+1}\right),
\ \
a\in\mathbb{R},
\
b>0.
\end{equation*}
\end{thm}

\begin{rem}
Here we provide some discussions about Theorem \ref{thmwsupscN4}.
\begin{enumerate}
[itemsep=0pt,
topsep=2pt,
parsep=0pt]

\item[(1)]
If $\alpha=0$ in \eqref{wsUPN4} (which reduces into the second order HyUP \eqref{sUP}), Theorem \ref{thmwsupscN4} turns into the case $N=4$ of \cite[Theorem 1.2]{Huang25}.

\item[(2)]
Notice that
\begin{align*}
&\left\{(N,\alpha):N\ge2,
\
\alpha>-1,
\
N\ge5\alpha+5\right\}\cap
\left\{(N,\alpha):N=4,
\
-1<\alpha\le0\right\}
\\&\quad=\left\{(N,\alpha):N=4,
\
-1<\alpha\le-\frac{1}{5}\right\},
\end{align*}
which means that the parameter range $(N,\alpha)$ of Theorem \ref{thmwsupsc-2} and that of Theorem \ref{thmwsupscN4} only intersect when $\left\{(N,\alpha):N=4,
\ -1<\alpha\le-\frac{1}{5}\right\}$.
\end{enumerate}
\end{rem}

The rest of this paper is organized as follows. In Section \ref{sect-2}, we aim to analyze the symmetry breaking phenomenon for the second order HyUP, and prove Theorem \ref{thm-3}. In Section \ref{sect-3}, we establish a family of weighted second order HyUP, and prove Theorems \ref{thmwsupsc-2} and \ref{thmwsupscN4}.

Throughout this paper, we will frequently use the following notations.
\begin{itemize}
\setlength{\itemsep}{0pt}
\setlength{\parsep}{0pt}
\setlength{\parskip}{2pt}

\item
$\mathbb{N}^*:=\{1,2,\cdots\}$ denotes the set which includes all natural numbers.

\item
$\mathbb{S}^{N-1}$ denotes the $(N\!-\!1)$-dimensional sphere with respect to the Hausdorff measure in $\mathbb{R}^N$.

\item $o_\varepsilon(1)$ denotes a quantity which tends to zero as $\varepsilon\to 0^+$.

\end{itemize}

\section{Preliminariles}

\noindent The spherical harmonics decomposition
discussed in Section \ref{sect-1.2} allows us to follow from
\cite[(2.3)$-$(2.5)]{Duong25-2} to obtain
the following identities.

\begin{lem}
[{\!\rm{\!\cite[(2.3)$-$(2.5)]{Duong25-2}}}]
\label{lemhduukd}
Let $N\geq 2$. For each
$u\in C^\infty_c(\mathbb{R}^N)$, $v_k$ is defined by \eqref{1031-2}:
\begin{align*}
\int_{\mathbb{R}^N}
\frac{|\Delta u|^2}{|x|^{2\alpha}}\mathrm{d}x
&=\sum^\infty_{k=0}
\left[\int_0^\infty |v_k''|^2r^{N+2k-2\alpha-1}
\mathrm{d}r
+(2\alpha+1)(N+2k-1)\int_0^\infty |v_k'|^2
r^{N+2k-2\alpha-3}
\mathrm{d}r\right],
\\
\int_{\mathbb{R}^N}\left|\nabla u\right|^2
\mathrm{d}x
&=\sum^\infty_{k=0}\int_0^\infty |v_k'|^2
r^{N+2k-1}\mathrm{d}r,
\\
\int_{\mathbb{R}^N}
\frac{\left|\nabla u\right|^2}{|x|^{\alpha+1}}
\mathrm{d}x
&=\sum^\infty_{k=0}
\left[\int_0^\infty |v_k'|^2r^{N+2k-\alpha-2}
\mathrm{d}r
+(\alpha+1)k\int_0^\infty |v_k|^2r^{N+2k-\alpha-4}
\mathrm{d}r
\right].
\end{align*}
\end{lem}

The Hardy inequality appearing in the following will be crucial to prove our main results.

\begin{lem}
[{\!\rm{\!\cite[(2.4)]{Huang25}}}]
\label{lem-2.3}
For any $\theta>1$ and $w\in C^1_c([0,\infty))$, then
\begin{equation}\label{HI}
\int_0^\infty |w'|^2 s^{\theta}\mathrm{d}s
\geq \frac{(\theta-1)^2}{4}
\int_0^\infty |w|^2 s^{\theta-2} \mathrm{d}s.
\end{equation}
\end{lem}

The inequality \eqref{HI} is the well-known weighted one-dimensional Hardy inequality. Here we also present another family of Hardy inequality.

\begin{lem}
[{\!\rm{\!\cite[Lemma 2.3]{Huang25}}}]
\label{lem-2.4}
For any $\theta\geq 1$, $\kappa\in\mathbb{R}$ and $w\in C^1_c([0,\infty))$, there holds
\begin{align}\label{HIk1}
\int_0^\infty |w'|^2 e^{-2s}s^{\theta} \mathrm{d}s
\geq -\theta \kappa\int_0^\infty |w|^2 e^{-2s} s^{\theta-1} \mathrm{d}s
-\kappa(\kappa-2)\int_0^\infty |w|^2 e^{-2s} s^{\theta} \mathrm{d}s.
\end{align}
\end{lem}

\section{Symmetry breaking phenomenon for the second order HyUP: proof of Theorem \ref{thm-3}}\label{sect-2}

\noindent In this section, we will show a symmetry breaking phenomenon for the second order HyUP \eqref{sUP}, that is, $C(N)<\frac{(N+1)^2}{4}$ for $N\in\{2,3\}$. Notice that the symmetry result heavily depends on the following conclusion obtained by Cazacu, Flynn and Lam  \cite[Lemma 3.4]{Cazacu22}: for any $N\geq 5$, \begin{equation*}
\inf_{k\in\mathbb{N}^*\cup\{0\}}
\frac{1}{\left[1+\frac{4k}
{(N+2k-3)^2}\right]^2}
\frac{(N+2k+1)^2}{4}
=\frac{(N+1)^2}{4}.
\end{equation*}
For $2\le N\le4$, we present the following calculation.

\begin{lem}\label{lemtr}
For $N\in\{2,3\}$, it holds that
\begin{equation*}
\inf_{k\in\mathbb{N}^*\cup\{0\}}
\frac{1}{\left[1+\frac{4k}
{(N+2k-3)^2}\right]^2}
\frac{(N+2k+1)^2}{4}
=\frac{(N-1)^4(N+3)^2}
{4(N^2-2N+5)^2}.
\end{equation*}
\end{lem}

\begin{proof}[\rm\textbf{Proof}]
Note that
\[
\frac{(N+2k+1)^2}{4\left[1+\frac{4k}
{(N+2k-3)^2}\right]^2}
=\frac{(N+2k-3)^4(N+2k+1)^2}
{4[(N+2k-3)^2+4k]^2}
:=\mathcal{J}(N,k).
\]

$\bullet$ \emph{Step 1: For $N\in\{2,3\}$, we will prove that}
\[
\inf_{k\in\mathbb{N}^*}
\mathcal{J}(N,k)
=\mathcal{J}(N,1).
\]

For $N\in\{2,3\}$, let us define the following function on the real axis $[2,\infty)$:
\begin{align*}
f(x)
:=\frac{(N+x-3)^4(N+x+1)^2}
{4[(N+x-3)^2+2x]^2}
=g(x)h(x),
\end{align*}
where
\begin{align*}
g(x)&:=\frac{(N+x+1)^2}
{4[x^2+2(N-2)x+(N-3)^2
]^{\frac{1}{4}}}>0,
\\
h(x)&:=\frac{(N+x-3)^4}
{[x^2+2(N-2)x+(N-3)^2
]^{\frac{7}{4}}}>0.
\end{align*}
It is easy to check that
\begin{align*}
g'(x)&
=\frac{(N+x+1)
\left[3x^2+3(2N-5)x+4(N-3)^2
-(N^2-N-2)\right]}
{8[x^2+2(N-2)x+(N-3)^2
]^{\frac{5}{4}}},
\\
h'(x)&=\frac{(N+x-3)^3
\left[x^2+(2N+3)x+(N^2-13N+30)
\right]}
{2[x^2+2(N-2)x+(N-3)^2
]^{\frac{11}{4}}}.
\end{align*}
Under the condition $N\in\{2,3\}$, it is easy to get that $-\frac{2N-5}{2}
\le\frac{1}{2}<2$, and then
\begin{align*}
&3x^2+3(2N-5)x+4(N-3)^2-(N^2-N-2)
\\&\quad\ge
12+6(2N-5)+4(N-3)^2-(N^2-N-2)
\\&\quad=
3N^2-11N+20>0,
\end{align*}
for all $x\ge2$. This leads to, $g'(x)>0$ for all $x\geq 2$.

On the other hand, using $N\in\{2,3\}$ again,  $-\frac{2N+3}{2}<0<2$ obviously, and then
\begin{align*}
&x^2+(2N+3)x+(N^2-13N+30)
\\&\quad\ge
4+2(2N+3)+(N^2-13N+30)
\\&\quad=N^2-9N+40>0,
\end{align*}
for all $x\ge2$. This implies that $h'(x)>0$ for all $x\geq 2$.

To sum up,
\[
f'(x)=g'(x)h(x)+g(x)h'(x)>0
\]
for all $x\geq 2$, this means that $f(x)$ is increasing for all $x\geq 2$. Thus, $f(x)\ge f(2)$ for all $x\ge2$, then for all $k\in\mathbb{N}^*$,
\begin{equation*}
\mathcal{J}(N,k)
=f(2k)\ge f(2)
=\mathcal{J}(N,1),
\end{equation*}
namely,
\begin{equation*}
\inf_{k\in \mathbb{N}^*}
\mathcal{J}(N,k)
=\mathcal{J}(N,1).
\end{equation*}
The proof of \emph{Step 1} finishes.

From \emph{Step 1}, we see that
\begin{align}\label{8.16-1}
\inf_{k\in\mathbb{N}^*\cup\{0\}}
\mathcal{J}(N,k)
&=\min\{\mathcal{J}(N,0),
\mathcal{J}(N,1)\}.
\end{align}

$\bullet$ \emph{Step 2: For $N\in\{2,3\}$, we will verify that}
\begin{equation}\label{8.16-2}
\mathcal{J}(N,1)
<\mathcal{J}(N,0).
\end{equation}

Notice that
\begin{eqnarray*}
\mathcal{J}(N,1)
=\frac{(N-1)^4(N+3)^2}
{4[(N-1)^2+4]^2}
=\left\{\arraycolsep=1.5pt
\begin{array}{ll}
\frac{1}{4},
\quad &\mbox{if}\  N=2,\\[2mm]
\frac{9}{4},
\quad &\mbox{if}\  N=3,
\end{array}
\right.
\end{eqnarray*}
and
\begin{eqnarray*}
\mathcal{J}(N,0)
=\frac{(N+1)^2}{4}
=\left\{\arraycolsep=1.5pt
\begin{array}{ll}
\frac{9}{4},
\quad &\mbox{if}\  N=2,\\[2mm]
4,
\quad &\mbox{if}\  N=3,
\end{array}
\right.
\end{eqnarray*}
then \eqref{8.16-2} holds. This ends the proof of \emph{Step 2}.

$\bullet$ \emph{Step 3: Conclusion.}

From \eqref{8.16-1} and \eqref{8.16-2},
\[
\inf_{k\in\mathbb{N}^*\cup\{0\}}
\mathcal{J}(N,k)
=\mathcal{J}(N,1)
=\frac{(N-1)^4(N+3)^2}
{4[(N-1)^2+4]^2},
\]
as our desired estimate. The proof is completed.
\end{proof}

\begin{rem}
\rm
Recall the definition
\begin{align*}
&\frac{(N+2k+1)^2}{4\left[1+\frac{4k}
{(N+2k-3)^2}\right]^2}
:=\mathcal{J}(N,k).
\end{align*}
From \cite[Lemma 3.4]{Cazacu22}, for all $N\ge5$,
\[
\inf_{k\in\mathbb{N}^*\cup\{0\}}
\mathcal{J}(N,k)
=\mathcal{J}(N,0).
\]
Meanwhile, from Lemma \ref{lemtr}, for $N\in\{2,3\}$,
\[
\inf_{k\in\mathbb{N}^*\cup\{0\}}
\mathcal{J}(N,k)
=\mathcal{J}(N,1).
\]
The above findings inspire us to consider the conjecture presented in \cite[Remark 2.4]{Cazacu22} may be not hold when $N\in\{2,3\}$.
\end{rem}

Now, based on the above conclusion, we are in a position to prove Theorem \ref{thm-3}.

\begin{proof}
[\rm\textbf{Proof of Theorem \ref{thm-3}}]
By the Hardy inequalities \eqref{HI} and the second-order HyUP on $\mathbb{R}^{N+2k}$, we have
\begin{align*}
C(N,k)
&=\inf_{u\in
W^{2,2}_{\mathrm{rad}}
(\mathbb{R}^{N+2k})
\setminus\{0\}}
\frac{\int_{\mathbb{R}^{N+2k}}
|\Delta u|^{2}\mathrm{d}x
\int_{\mathbb{R}^{N+2k}}
|\nabla u|^{2}\mathrm{d}x}
{\left[\int_{\mathbb{R}^{N+2k}}
\left(\frac{|\nabla u|^{2}}{|x|}
+k\frac{|u|^{2}}{|x|^3}\right)
\mathrm{d}x\right]^2}
\nonumber\\&\ge\frac{1}{\left[1+\frac{4k}
{(N+2k-3)^2}\right]^2}
\inf_{u\in
W^{2,2}_{\mathrm{rad}}
(\mathbb{R}^{N+2k})
\setminus\{0\}}
\frac{\int_{\mathbb{R}^{N+2k}}
|\Delta u|^{2}\mathrm{d}x
\int_{\mathbb{R}^{N+2k}}
|\nabla u|^{2}\mathrm{d}x}
{\left(\int_{\mathbb{R}^{N+2k}}
\frac{|\nabla u|^{2}}{|x|}
\mathrm{d}x\right)^2}
\nonumber\\&=
\frac{(N+2k+1)^2}{4\left[1+\frac{4k}
{(N+2k-3)^2}\right]^2},
\end{align*}
for all $k\in\mathbb{N}^*$. Thus, from Lemma \ref{lemtr},
\begin{equation*}
C(N)
=\inf_{k\in\mathbb{N}^*\cup\{0\}}
C(N,k)
\ge\inf_{k\in\mathbb{N}^*\cup\{0\}}
\frac{(N+2k+1)^2}{4\left[1+\frac{4k}
{(N+2k-3)^2}\right]^2}
=\frac{(N-1)^4(N+3)^2}
{4[(N-1)^2+4]^2}.
\end{equation*}

Next, by testing the constant $C(N,1)$ with $U(x)=\exp(-|x|)$, we deduce that
\begin{align*}
\int_{\mathbb{R}^{N+2}}
|\Delta U|^2\mathrm{d}x
&=\left|\mathbb{S}^{N+1}\right|
\int_0^\infty\left|U''
+\frac{N+1}{r}U'\right|^2 r^{N+1}\mathrm{d}r
\\&=\left|\mathbb{S}^{N+1}\right|
\left[\int_0^\infty|U''|^2 r^{N+1}\mathrm{d}r
+(N+1)\int_0^\infty|U'|^2 r^{N-1}\mathrm{d}r\right]
\\&=\left|\mathbb{S}^{N+1}\right|
\left[\int_0^\infty e^{-2r} r^{N+1}\mathrm{d}r
+(N+1)\int_0^\infty e^{-2r} r^{N-1}\mathrm{d}r\right]
\\&=\left|\mathbb{S}^{N+1}\right|
\left[\frac{\Gamma(N+2)}{2^{N+2}}
+\frac{(N+1)\Gamma(N)}{2^{N}}
\right],
\end{align*}
and
\begin{align*}
\int_{\mathbb{R}^{N+2}}
\left|\nabla U\right|^2 \mathrm{d}x
&=\left|\mathbb{S}^{N+1}\right|
\int_0^\infty|U'|^2r^{N+1}\mathrm{d}r
=\left|\mathbb{S}^{N+1}\right|
\int_0^\infty e^{-2r} r^{N+1}\mathrm{d}r
=\left|\mathbb{S}^{N+1}\right|
\frac{\Gamma(N+2)}{2^{N+2}},\\
\int_{\mathbb{R}^{N+2}}
\frac{\left|\nabla U\right|^2}{|x|} \mathrm{d}x
&=\left|\mathbb{S}^{N+1}\right|
\int_0^\infty|U'|^2r^{N}\mathrm{d}r
=\left|\mathbb{S}^{N+1}\right|
\int_0^\infty e^{-2r} r^{N}\mathrm{d}r
=\left|\mathbb{S}^{N+1}\right|
\frac{\Gamma(N+1)}{2^{N+1}},\\
\int_{\mathbb{R}^{N+2}}
\frac{\left|U\right|^2}
{|x|^3}\mathrm{d}x
&=\left|\mathbb{S}^{N+1}\right|
\int_0^\infty|U|^2 r^{N-2}\mathrm{d}r
=\left|\mathbb{S}^{N+1}\right|
\int_0^\infty e^{-2r} r^{N-2}\mathrm{d}r
=\left|\mathbb{S}^{N+1}\right|
\frac{\Gamma(N-1)}{2^{N-1}}.
\end{align*}
Therefore,
\begin{align*}
\frac{\int_{\mathbb{R}^{N+2}}
|\Delta U|^2\mathrm{d}x
\int_{\mathbb{R}^{N+2}}
\left|\nabla U\right|^2 \mathrm{d}x}{
\left(\int_{\mathbb{R}^{N+2}}
\frac{\left|\nabla U\right|^2}
{|x|}+\frac{\left|U\right|^2}
{|x|^3}\mathrm{d}x\right)^2}
&=\frac{
\left[\frac{\Gamma(N+2)}{2^{N+2}}
+\frac{(N+1)\Gamma(N)}{2^{N}}
\right]
\frac{\Gamma(N+2)}{2^{N+2}}}
{\left[\frac{\Gamma(N+1)}{2^{N+1}}
+\frac{\Gamma(N-1)}{2^{N-1}}
\right]^2}
=\frac{N(N+4)(N^2-1)^2}
{4(N^2-N+4)^2},
\end{align*}
with the help of the property of the Gamma function: $\Gamma(t+1)=t\Gamma(t)$ for $t>0$.
Thus,
\begin{align*}
C(N,1)
&\le\frac{\int_{\mathbb{R}^{N+2}}
|\Delta U|^{2}\mathrm{d}x
\int_{\mathbb{R}^{N+2}}
|\nabla U|^{2}\mathrm{d}x}
{\left[\int_{\mathbb{R}^{N+2}}
\left(\frac{|\nabla U|^{2}}{|x|}
+\frac{|U|^{2}}{|x|^3}\right)
\mathrm{d}x\right]^2}
\nonumber\\&=
\frac{N(N+4)(N^2-1)^2}
{4(N^2-N+4)^2}
\nonumber\\&<\frac{(N+1)^4(N+5)^2}
{4[(N+1)^2+8]^2}
\le\frac{(N+2k+1)^2}{4\left[1+\frac{4k}
{(N+2k-3)^2}\right]^2}
\le C(N,k),
\ \
\mathrm{for}
\
\mathrm{all}
\
k\ge2,
\end{align*}
where the third inequality follows from the fact that the function
\[
\frac{(N+2k+1)^2}{4\left[1+\frac{4k}
{(N+2k-3)^2}\right]^2}
\]
is positive and increasing for all $k\in\mathbb{N}^*$. As a consequence of the above inequality,
\begin{equation*}
C(N,1)
\le
\frac{N(N+4)(N^2-1)^2}
{4(N^2-N+4)^2}
<\frac{(N+1)^2}{4}=C(N,0).
\end{equation*}
Therefore, $C(N)=C(N,1)$, and then
\[
C(N)
\le
\frac{N(N+4)(N^2-1)^2}
{4(N^2-N+4)^2}
<\frac{(N+1)^2}{4}.
\]
The proof is completed here.
\end{proof}

\section{Symmetry phenomenon for the weighted second order HyUP: proof of Theorems \ref{thmwsupsc-2} and \ref{thmwsupscN4}}
\label{sect-3}

\noindent In this section, we consider the sharp weighted second order HyUP, and prove Theorems \ref{thmwsupsc-2} and \ref{thmwsupscN4}.

\subsection{The case $N\ge2$, $\alpha>-1$ and $N\ge5\alpha+5$: proof of Theorem \ref{thmwsupsc-2}}
\label{sect-3.2}

\noindent In this subsection, we prove Theorem \ref{thmwsupsc-2}. Now, inspired by the argument
of Cazacu, Flynn and Lam in \cite[Lemma 3.4]{Cazacu22}, we establish the following lemma which will be applied to compute the sharp constant of \eqref{wsUP}.

\begin{lem}\label{lemtrw}
Let $N\geq 2$ and $\alpha>-1$ satisfy
$N\geq 5\alpha+5$, then
\begin{equation*}
\inf_{k\in\mathbb{N}^*\cup\{0\}}
\frac{1}
{\left[1+\frac{4(\alpha+1)k}
{(N+2k-\alpha-3)^2}\right]^2}
\left(\frac{N+2k+3\alpha+1}{2}
\right)^2
=\frac{(N+3\alpha+1)^2}{4}.
\end{equation*}
\end{lem}

\begin{proof}[\rm\textbf{Proof}]
Note that
\begin{align*}
&\frac{1}
{\left[1+\frac{4(\alpha+1)k}
{(N+2k-\alpha-3)^2}\right]^2}
\left(\frac{N+2k+3\alpha+1}{2}
\right)^2
\\&\quad
=\frac{(N+2k-\alpha-3)^4
(N+2k+3\alpha+1)^2}
{4\left[(N+2k-\alpha-3)^2
+4(\alpha+1)k\right]^2}
:=\mathcal{K}(N,\alpha,k).
\end{align*}

$\bullet$ \emph{Step 1: For $N\ge2$
and $\alpha>-1$ satisfying $N\ge 5\alpha+5$,
we will verify that}
\begin{equation*}
\inf_{k\in\mathbb{N}^*}
\mathcal{K}(N,\alpha,k)
=\mathcal{K}(N,\alpha,1).
\end{equation*}

Under the assumption $N\ge2$ and $\alpha>-1$
satisfying $N\ge5\alpha+5$, we first show that
\begin{equation*}
\mathcal{F}(x)
:=\frac{(N+x-\alpha-3)^4
(N+x+3\alpha+1)^2}
{4\left[(N+x-\alpha-3)^2
+2(\alpha+1)x\right]^2}
\end{equation*}
is non-decreasing for $x\ge2$. In fact,
let $t:=N+x-\alpha-3\ge 4\alpha+4>0$
(due to $N\ge5\alpha+5$, $x\ge2$ and $\alpha>-1$),
consider the function $\mathcal{G}:[2,\infty)
\times[4\alpha+4,\infty)$:
\begin{equation*}
\mathcal{G}(x,t)
=\frac{t^4(t+4\alpha+4)^2}
{4\left[t^2+2(\alpha+1)x\right]^2}.
\end{equation*}
The direct computation yields
\begin{align*}
\frac{\partial\mathcal{G}}
{\partial x}(x,t)
&=-\frac{(\alpha+1)t^4
(t+4\alpha+4)^2}
{\left[t^2+2(\alpha+1)x\right]^3},
\\\frac{\partial\mathcal{G}}
{\partial t}(x,t)
&=\frac{(t+4\alpha+4)t^3
\left[t^3+2(\alpha+1)x
(3t+8\alpha+8)\right]}
{2\left[t^2+2(\alpha+1)x\right]^3}.
\end{align*}
Note that $\mathcal{F}(x)=\mathcal{G}(x,t)$
with $t=N+x-\alpha-3$, then we get, for $x\geq 2$
and $t\geq 4\alpha+4>0$,
\begin{align*}
\mathcal{F}'(x)&=
\frac{\partial\mathcal{G}}
{\partial x}(x,t)
+\frac{\partial\mathcal{G}}
{\partial t}(x,t)
\\& =-\frac{(\alpha+1)t^4(t+4\alpha+4)^2}
{\left[t^2+2(\alpha+1)x\right]^3}
+\frac{(t+4\alpha+4)t^3
\left[t^3+2(\alpha+1)x
(3t+8\alpha+8)\right]}
{2\left[t^2
+2(\alpha+1)x\right]^3}
\\&=\frac{(t+4\alpha+4)t^3
\left[t^3+2(\alpha+1)x
(3t+8\alpha+8)\right]
-2(\alpha+1)t^4(t+4\alpha+4)^2}
{2\left[t^2+2(\alpha+1)x\right]^3}
\\&\geq
\frac{(t+4\alpha+4)t^6
-2(\alpha+1)t^4(t+4\alpha+4)^2}
{2\left[t^2+2(\alpha+1)x\right]^3}
\\&=\frac{(t+4\alpha+4)t^4
[t^2-2(\alpha+1)(t+4\alpha+4)]}
{2\left[t^2+2(\alpha+1)x\right]^3}
\\&=\frac{(t+4\alpha+4)t^4
[t+2(\alpha+1)][t-4(\alpha+1)]}
{2\left[t^2+2(\alpha+1)x\right]^3}
\\&\geq 0,
\end{align*}
which means that $\mathcal{F}(x)$ is
non-decreasing for $x\ge2$, and then $\mathcal{F}(x)
\ge\mathcal{F}(2)$ for all $x\ge2$.
Thus, for all $k\in\mathbb{N}^*$,
\begin{equation}
\mathcal{K}(N,\alpha,k)
=\mathcal{F}(2k)
\ge\mathcal{F}(2)
=\mathcal{K}(N,\alpha,1),
\end{equation}
that is,
\begin{equation*}
\inf_{k\in\mathbb{N}^*}
\mathcal{K}(N,\alpha,k)
=\mathcal{K}(N,\alpha,1).
\end{equation*}
In view of this, the proof of \emph{Step 1}
is completed.

From \emph{Step 1},
\begin{equation}\label{8.11-1}
\inf_{k\in\mathbb{N}^*\cup\{0\}}
\mathcal{K}(N,\alpha,k)
=\min\left\{\mathcal{K}(N,\alpha,0),
\mathcal{K}(N,\alpha,1)\right\}.
\end{equation}

$\bullet$ \emph{Step 2: For $N\geq 2$
and $\alpha>-1$ satisfying $N\geq 5\alpha+5$,
we will show that}
\begin{equation*}
\mathcal{K}(N,\alpha,1)
\ge\mathcal{K}(N,\alpha,0).
\end{equation*}

The above inequality is equivalent to
\begin{equation}\label{8.11-2}
\frac{(N-\alpha-1)^4(N+3\alpha+3)^2}
{4\left[(N-\alpha-1)^2
+4(\alpha+1)\right]^2}
\ge\frac{(N+3\alpha+1)^2}{4}.
\end{equation}
Let $y:=N-\alpha-1\geq 4\alpha+4>0$
(due to $N\ge5\alpha+5$ and $\alpha>-1$),
then we can rewrite \eqref{8.11-2} as
\begin{equation}\label{8.16-3}
y^4(y+4\alpha+4)^2
\geq (y+4\alpha+2)^2(y^2+4\alpha+4)^2.
\end{equation}
Using $N\ge2$, $\alpha>-1$ and $y=N-\alpha-1>0$
again, we see that
\begin{align*}
\begin{cases}
y^2>0,\\
y+4\alpha+4=N+3\alpha+3>0,\\
y+4\alpha+2=N+3\alpha+1>0,\\
y^2+4\alpha+4>0,
\end{cases}
\end{align*}
then \eqref{8.16-3} reduces into
\begin{equation}\label{8.10-1}
y^2\ge2(\alpha+1)y
+4(\alpha+1)(2\alpha+1).
\end{equation}
In order to complete the proof of \emph{Step 2},
it remains to verify that \eqref{8.10-1} holds.

Indeed, by $y\geq 4\alpha+4$ and $\alpha>-1$, then
\begin{equation*}
y^2-2(\alpha+1)y
=[y-(\alpha+1)]^2-(\alpha+1)^2
\geq 8(\alpha+1)^2\geq 4(\alpha+1)(2\alpha+1),
\end{equation*}
which means that \eqref{8.10-1} holds.
This finishes the proof of \emph{Step 2}.

$\bullet$ \emph{Step 3: Conclusion.}

Combining \eqref{8.11-1} and \emph{Step 2},
\begin{equation*}
\inf_{k\in\mathbb{N}^*\cup\{0\}}
\mathcal{K}(N,\alpha,k)
=\mathcal{K}(N,\alpha,0)
=\frac{(N+3\alpha+1)^2}{4},
\end{equation*}
as our desired estimate.
\end{proof}

Now, based on the above preparations,
we follow the approach of Cazacu, Flynn
and Lam \cite[Proof of Theorem 2.3]{Cazacu22}
to prove Theorem \ref{thmwsupsc-2}.

\begin{proof}
[\rm\textbf{Proof of Theorem \ref{thmwsupsc-2}}]

The proof follows in two steps as follows.

$\bullet$ \emph{Step 1: We will prove
that the weighted second order HyUP
\eqref{wsUP} holds}.

Let $u\in C^\infty_c(\mathbb{R}^N)$,
using Lemma \ref{lemhduukd},
the inequality \eqref{wsUP} is equivalent to
\begin{align}\label{8.11-3}
&\sum^\infty_{k=0}
\left[\int_0^\infty |v_k''|^2r^{N+2k-2\alpha-1}
\mathrm{d}r
+(2\alpha+1)(N+2k-1)
\int_0^\infty |v_k'|^2r^{N+2k-2\alpha-3}
\mathrm{d}r
\right]
\nonumber\\
&\qquad \times\sum^\infty_{k=0}
\int_0^\infty |v_k'|^2r^{N+2k-1}\mathrm{d}r
\nonumber\\
&\quad\ge C(N,\alpha)
\left\{\sum^\infty_{k=0}
\left[\int_0^\infty |v_k'|^2r^{N+2k-\alpha-2}
\mathrm{d}r
+(\alpha+1)k\int_0^\infty |v_k|^2r^{N+2k-\alpha-4}
\mathrm{d}r
\right]\right\}^2.
\end{align}
By Cauchy-Schwarz inequality, in order to
verify \eqref{8.11-3}, it suffices to show that,
for all $k\in\mathbb{N}^*\cup\{0\}$,
\begin{align}\label{8.11-4}
&\left[\int_0^\infty |v_k''|^2r^{N+2k-2\alpha-1}
\mathrm{d}r
+(2\alpha+1)(N+2k-1)\int_0^\infty |v_k'|^2
r^{N+2k-2\alpha-3}
\mathrm{d}r
\right]\int_0^\infty |v_k'|^2r^{N+2k-1}\mathrm{d}r
\nonumber\\
&\quad \ge C(N,\alpha)
\left[\int_0^\infty |v_k'|^2r^{N+2k-\alpha-2}
\mathrm{d}r
+(\alpha+1)k\int_0^\infty |v_k|^2r^{N+2k-\alpha-4}
\mathrm{d}r
\right]^2.
\end{align}
For $k\in\mathbb{N}^*$, choosing $\theta=N+2k-\alpha-2$ in the Hardy inequality \eqref{HI}, we have
\begin{equation*}
\int_0^\infty |v_k|^2r^{N+2k-\alpha-4}
\mathrm{d}r
\le\frac{4}{(N+2k-\alpha-3)^2}
\int_0^\infty |v_k'|^2r^{N+2k-\alpha-2}
\mathrm{d}r,
\end{equation*}
this gives
\begin{align*}
&\int_0^\infty |v_k'|^2r^{N+2k-\alpha-2}
\mathrm{d}r
+(\alpha+1)k\int_0^\infty |v_k|^2r^{N+2k-\alpha-4}
\mathrm{d}r
\\&\quad\le
\left[1+\frac{4(\alpha+1)k}
{(N+2k-\alpha-3)^2}\right]
\int_0^\infty |v_k'|^2r^{N+2k-\alpha-2}
\mathrm{d}r.
\end{align*}
Thus, to verify \eqref{8.11-4},
it is enough to prove that,
for any $k\in\mathbb{N}^*\cup\{0\}$,
\begin{align*}
&\left[\int_0^\infty |v_k''|^2r^{N+2k-2\alpha-1}
\mathrm{d}r
+(2\alpha+1)(N+2k-1)\int_0^\infty |v_k'|^2
r^{N+2k-2\alpha-3}
\mathrm{d}r
\right]\int_0^\infty |v_k'|^2r^{N+2k-1}\mathrm{d}r
\nonumber\\
&\quad \ge C(N,\alpha)
\left[1+\frac{4(\alpha+1)k}
{(N+2k-\alpha-3)^2}\right]^2
\left(\int_0^\infty |v_k'|^2r^{N+2k-\alpha-2}
\mathrm{d}r\right)^2.
\end{align*}
In view of this, we try to get
the sharp constant $D(N,\alpha,k)$
in the following inequality:
\begin{align}\label{1dni}
&
\left[\int_0^\infty |v_k''|^2r^{N+2k-2\alpha-1}
\mathrm{d}r
+(2\alpha+1)(N+2k-1)\int_0^\infty |v_k'|^2
r^{N+2k-2\alpha-3}
\mathrm{d}r
\right]\int_0^\infty |v_k'|^2r^{N+2k-1}\mathrm{d}r
\nonumber\\
&\quad \ge D(N,\alpha,k)
\left(\int_0^\infty|v_k'|^2
r^{N+2k-\alpha-2}
\mathrm{d}r\right)^2,
\end{align}
for any $k\in\mathbb{N}^*\cup\{0\}$. Then
\begin{equation*}
C(N,\alpha)
\geq\inf_{k\in\mathbb{N}^*\cup\{0\}}
\frac{D(N,\alpha,k)}
{\left[1+\frac{4(\alpha+1)k}
{(N+2k-\alpha-3)^2}\right]^2}.
\end{equation*}
Let $w:=v_k'$, then \eqref{1dni} becomes,
for any $k\in\mathbb{N}^*\cup\{0\}$,
\begin{align*}
&\left[\int_0^\infty |w'|^2r^{N+2k-2\alpha-1}
\mathrm{d}r
+(2\alpha+1)(N+2k-1)\int_0^\infty |w|^2
r^{N+2k-2\alpha-3}\mathrm{d}r
\right]\int_0^\infty |w|^2r^{N+2k-1}\mathrm{d}r
\nonumber\\
&\quad \ge D(N,\alpha,k)
\left(\int_0^\infty|w|^2r^{N+2k-\alpha-2}
\mathrm{d}r\right)^2.
\end{align*}
Fix $0<\varepsilon\ll 1$.
Let $w(r):=r^{2\alpha+1}z(r)$ for
$r\in[\varepsilon,\infty)$.
Then applying integration by parts, we see that
\begin{align*}
&\int_\varepsilon^\infty |w'|^2r^{N+2k-2\alpha-1}
\mathrm{d}r
+(2\alpha+1)(N+2k-1)
\int_\varepsilon^\infty |w|^2r^{N+2k-2\alpha-3}
\mathrm{d}r
\\&\quad=
\int_\varepsilon^\infty \left|(2\alpha+1)r^{2\alpha}z
+r^{2\alpha+1}z'\right|^2
r^{N+2k-2\alpha-1}
\mathrm{d}r
\nonumber\\&\qquad
+(2\alpha+1)(N+2k-1)
\int_\varepsilon^\infty |z|^2r^{N+2k+2\alpha-1}
\mathrm{d}r
\\&\quad=
(2\alpha+1)^2
\int_\varepsilon^\infty \left|z\right|^2
r^{N+2k+2\alpha-1}
\mathrm{d}r
+\int_\varepsilon^\infty
\left|z'\right|^2
r^{N+2k+2\alpha+1}
\mathrm{d}r
\nonumber\\&\qquad
+2(2\alpha+1)
\int_\varepsilon^\infty
zz'r^{N+2k+2\alpha}
\mathrm{d}r
+(2\alpha+1)(N+2k-1)
\int_\varepsilon^\infty |z|^2r^{N+2k+2\alpha-1}
\mathrm{d}r
\\&\quad=
(2\alpha+1)^2
\int_\varepsilon^\infty \left|z\right|^2
r^{N+2k+2\alpha-1}
\mathrm{d}r
+\int_\varepsilon^\infty
\left|z'\right|^2
r^{N+2k+2\alpha+1}
\mathrm{d}r
\nonumber\\&\qquad
-(2\alpha+1)(N+2k+2\alpha)
\int_\varepsilon^\infty
|z|^2
r^{N+2k+2\alpha-1}
\mathrm{d}r+o_\varepsilon(1)
\nonumber\\&\qquad
+(2\alpha+1)(N+2k-1)
\int_\varepsilon^\infty |z|^2r^{N+2k+2\alpha-1}
\mathrm{d}r
\\&\quad=
\int_\varepsilon^\infty |z'|^2r^{N+2k+2\alpha+1}
\mathrm{d}r
+o_\varepsilon(1).
\end{align*}
Hence, using integrating by parts and
H\"{o}lder inequality, we get
\begin{align}\label{1dnibk}
&
\left[\int_\varepsilon^\infty |w'|^2
r^{N+2k-2\alpha-1}
\mathrm{d}r
+(2\alpha+1)(N+2k-1)
\int_\varepsilon^\infty |w|^2r^{N+2k-2\alpha-3}
\mathrm{d}r
\right]\int_\varepsilon^\infty |w|^2
r^{N+2k-1}\mathrm{d}r
\nonumber\\
&\quad =\left[\int_\varepsilon^\infty |z'|^2
r^{N+2k+2\alpha+1}
\mathrm{d}r
+o_\varepsilon(1)\right]
\int_\varepsilon^\infty |z|^2r^{N+2k+4\alpha+1}
\mathrm{d}r
\nonumber\\
&\quad =\int_\varepsilon^\infty |z'|^2
r^{N+2k+2\alpha+1}
\mathrm{d}r
\int_\varepsilon^\infty |z|^2r^{N+2k+4\alpha+1}
\mathrm{d}r
+o_\varepsilon(1)
\nonumber\\
&\quad \ge \left(\int_\varepsilon^\infty zz'
r^{N+2k+3\alpha+1}
\mathrm{d}r\right)^2
+o_\varepsilon(1)
\nonumber\\
&\quad =\frac{1}{4}
\left[\int_\varepsilon^\infty
\left(|z(r)|^2\right)'
r^{N+2k+3\alpha+1}
\mathrm{d}r\right]^2
+o_\varepsilon(1)
\nonumber\\
&\quad =\frac{(N+2k+3\alpha+1)^2}{4}
\left(\int_\varepsilon^\infty
|z(r)|^2 r^{N+2k+3\alpha}
\mathrm{d}r\right)^2
+o_\varepsilon(1)
\nonumber\\
&\quad =\frac{(N+2k+3\alpha+1)^2}{4}
\left(\int_\varepsilon^\infty |w|^2r^{N+2k-\alpha-2}
\mathrm{d}r\right)^2
+o_\varepsilon(1).
\end{align}
Let $\varepsilon\to 0$, we obtain,
for any $k\in\mathbb{N}^*\cup\{0\}$,
\begin{align*}
&\left[\int_0^\infty |w'|^2r^{N+2k-2\alpha-1}
\mathrm{d}r
+(2\alpha+1)(N+2k-1)\int_0^\infty |w|^2
r^{N+2k-2\alpha-3}
\mathrm{d}r
\right]\int_0^\infty |w|^2r^{N+2k-1}\mathrm{d}r
\nonumber\\
&\quad \ge \frac{(N+2k+3\alpha+1)^2}{4}
\left(\int_0^\infty |w|^2r^{N+2k-\alpha-2}
\mathrm{d}r\right)^2.
\end{align*}
Thus, under the assumptions
$N\geq 2$ and $\alpha>-1$
satisfying $N\geq 5\alpha+5$,
we deduce from Lemma \ref{lemtrw} that
\begin{align*}
C(N,\alpha)
&\geq\inf_{k\in\mathbb{N}^*\cup\{0\}}
\frac{1}
{\left[1+\frac{4(\alpha+1)k}
{(N+2k-\alpha-3)^2}\right]^2}
\frac{(N+2k+3\alpha+1)^2}{4}
\\&=\frac{(N+3\alpha+1)^2}{4}.
\end{align*}
Therefore, by standard density argument,
for all $u\in \mathcal{H}^2_{\alpha,0}(\mathbb{R}^N)$,
the desired result \eqref{wsUP} follows.

$\bullet$ \emph{Step 2: We will verify the attainability of the sharp constant
$\frac{(N+3\alpha+1)^2}
{4}$ of \eqref{wsUP}.}

From Lemma \ref{lemtrw}, we see that
\begin{align*}
&\inf_{k\in\mathbb{N}^*\cup\{0\}}
\frac{1}
{\left[1+\frac{4(\alpha+1)k}
{(N+2k-\alpha-3)^2}\right]^2}
\left(\frac{N+2k+3\alpha+1}{2}
\right)^2
\\&\quad
=\inf_{k\in\mathbb{N}^*\cup\{0\}}
\frac{1}
{\left[1+\frac{4(\alpha+1)k}
{(N+2k-\alpha-3)^2}\right]^2}
\left(\frac{N+2k+3\alpha+1}{2}
\right)^2\bigg|_{k=0}
\\&\quad=\frac{(N+3\alpha+1)^2}{4},
\end{align*}
then if the extremal function of \eqref{wsUP} exists,
it must be radial (according to $k=0$).
Moreover, if $u(|x|)$ is the extremal
function of \eqref{wsUP}, then $w(r):=u'(r)$
also must be the extremal function of
\begin{align}\label{wuppks}
&\left[
\int^\infty_0\left|w'\right|^2
r^{N-2\alpha-1}\mathrm{d}r
+(2\alpha+1)(N-1)
\int^\infty_0\left|w\right|^2
r^{N-2\alpha-3}\mathrm{d}r\right]
\int^\infty_0\left|w\right|^2
r^{N-1}\mathrm{d}r
\nonumber
\\
&\quad \geq \frac{(N+3\alpha+1)^2}{4}
\left(\int^\infty_0\left|w\right|^2
r^{N-\alpha-2}\mathrm{d}r\right)^2.
\end{align}
By \eqref{1dnibk} with $k=0$,
we know that the equality in \eqref{wuppks}
holds when applying H\"{o}lder inequality to
\eqref{1dnibk} if and only if
\begin{align*}
|z'|^2 r^{N+2\alpha+1}
=C_1|z|^2 r^{N+4\alpha+1}
&\Rightarrow
z(r)=C_2
\exp\left(-cr^{\alpha+1}\right)
\\&\Rightarrow
w(r)=C_2r^{2\alpha+1}
\exp\left(-cr^{\alpha+1}\right)
\ \ \ \ \ \ \ \ \ \,
(\mbox{by}
\
w(r)=r^{2\alpha+1}z(r))
\\&\Rightarrow
u'(r)=C_2r^{2\alpha+1}
\exp\left(-cr^{\alpha+1}\right)
\qquad\qquad\ \ \ \,
(\mbox{by}
\
w(r)=u'(r)),
\end{align*}
for some $c,C_1>0$ and $C_2\in\mathbb{R}$.
Therefore, the extremal function of \eqref{wsUP}
must be of the form $\bar{u}(x)=a(1+b|x|^{\alpha+1})
\exp\left(-b |x|^{\alpha+1}\right)$
with $a\in \mathbb{R}$ and $b>0$.

Indeed, direct computations yield
\begin{align*}
\int_{\mathbb{R}^N}\frac{|\Delta \bar{u}|^2}
{|x|^{2\alpha}}
\mathrm{d}x
&=\left|\mathbb{S}^{N-1}\right|
\int^\infty_0
\left|\bar{u}''+\frac{N-1}{r}
\bar{u}'\right|^2
r^{N-2\alpha-1}\mathrm{d}r
\\&=\left|\mathbb{S}^{N-1}\right|
a^2b^4(\alpha+1)^2
\int^\infty_0
\left|(N+2\alpha)-b(\alpha+1)
r^{\alpha+1}\right|^2
\exp\left(-2br^{\alpha+1}\right)
r^{N+2\alpha-1}\mathrm{d}r
\\&=\left|\mathbb{S}^{N-1}\right|
a^2b^4(\alpha+1)^2
\bigg[(N+2\alpha)^2
\int^\infty_0
\exp\left(-2br^{\alpha+1}\right)
r^{N+2\alpha-1}\mathrm{d}r
\\&\quad
-2(N+2\alpha)b(\alpha+1)
\int^\infty_0
\exp\left(-2br^{\alpha+1}\right)
r^{N+3\alpha}
\mathrm{d}r
\nonumber\\&\quad
+b^2(\alpha+1)^2
\int^\infty_0
\exp\left(-2br^{\alpha+1}\right)
r^{N+4\alpha+1}\mathrm{d}r\bigg].
\end{align*}
Next, we change the variable: $s=2br^{\alpha+1}$,
then
\begin{align*}
\int^\infty_0
\exp\left(-2br^{\alpha+1}\right)
r^{N+2\alpha-1}
\mathrm{d}r
&=\frac{1}{(\alpha+1)
(2b)^{\frac{N+2\alpha}{\alpha+1}}}
\int^\infty_0
e^{-s}
s^{\frac{N+2\alpha}{\alpha+1}-1}
\mathrm{d}s
\\&=\frac{1}{(\alpha+1)
(2b)^{\frac{N+2\alpha}{\alpha+1}}}
\Gamma\left(\frac{N+2\alpha}
{\alpha+1}\right),
\\
\int^\infty_0
\exp\left(-2br^{\alpha+1}\right)
r^{N+3\alpha}\mathrm{d}r
&=\frac{1}{(\alpha+1)
(2b)^{\frac{N+2\alpha}{\alpha+1}+1}}
\int^\infty_0
e^{-s}s^{\frac{N+2\alpha}{\alpha+1}}
\mathrm{d}s
\\&=\frac{1}{(\alpha+1)
(2b)^{\frac{N+2\alpha}{\alpha+1}+1}}
\Gamma\left(\frac{N+2\alpha}
{\alpha+1}+1\right),
\\
\int^\infty_0
\exp\left(-2br^{\alpha+1}\right)
r^{N+4\alpha+1}\mathrm{d}r
&=\frac{1}{(\alpha+1)
(2b)^{\frac{N+2\alpha}{\alpha+1}+2}}
\int^\infty_0
e^{-s}s^{\frac{N+2\alpha}{\alpha+1}+1}
\mathrm{d}s
\\&=\frac{1}{(\alpha+1)
(2b)^{\frac{N+2\alpha}{\alpha+1}+2}}
\Gamma\left(\frac{N+2\alpha}
{\alpha+1}+2\right).
\end{align*}
Hence, using the property of the Gamma function:
$\Gamma(t+1)=t\Gamma(t)$ for $t>0$, there holds
\begin{align*}
\int_{\mathbb{R}^N}
\frac{|\Delta \bar{u}|^2}{|x|^{2\alpha}}
\mathrm{d}x
=\frac{\left|\mathbb{S}^{N-1}\right|
a^2b^4(\alpha+1)}
{4(2b)^{\frac{N+2\alpha}{\alpha+1}}}
(N+3\alpha+1)(N+2\alpha)
\Gamma\left(\frac{N+2\alpha}
{\alpha+1}\right).
\end{align*}
Also,
\begin{align*}
\int_{\mathbb{R}^N}
|\nabla \bar{u}|^2\mathrm{d}x
&=\left|\mathbb{S}^{N-1}\right|
\int^\infty_0
\left|\bar{u}'\right|^2
r^{N-1}\mathrm{d}r
\\&=\left|\mathbb{S}^{N-1}\right|
a^2b^4(\alpha+1)^2
\int^\infty_0
\exp\left(-2br^{\alpha+1}\right)
r^{N+4\alpha+1}\mathrm{d}r
\\&=
\frac{\left|\mathbb{S}^{N-1}\right|
a^2b^4(\alpha+1)}
{(2b)^{\frac{N+2\alpha}{\alpha+1}+2}}
\Gamma\left(\frac{N+2\alpha}
{\alpha+1}+2\right),
\end{align*}
and
\begin{align*}
\int_{\mathbb{R}^N}
\frac{|\nabla \bar{u}|^2}{|x|^{\alpha+1}}
\mathrm{d}x
&=\left|\mathbb{S}^{N-1}\right|
\int^\infty_0
\left|\bar{u}'\right|^2
r^{N-\alpha-2}\mathrm{d}r
\\&=\left|\mathbb{S}^{N-1}\right|
a^2b^4(\alpha+1)^2
\int^\infty_0
\exp\left(-2br^{\alpha+1}\right)
r^{N+3\alpha}\mathrm{d}r
\\&=
\frac{\left|\mathbb{S}^{N-1}\right|
a^2b^4(\alpha+1)}
{(2b)^{\frac{N+2\alpha}{\alpha+1}+1}}
\Gamma\left(\frac{N+2\alpha}
{\alpha+1}+1\right).
\end{align*}
Therefore,
\begin{align*}
\dfrac{\int_{\mathbb{R}^N}
\frac{|\Delta \bar{u}|^2}{|x|^{2\alpha}}
\mathrm{d}x
\int_{\mathbb{R}^N}
\left|\nabla \bar{u}\right|^2 \mathrm{d}x}
{\left(\int_{\mathbb{R}^N}
\frac{\left|\nabla \bar{u}\right|^2}{|x|^{\alpha+1}}
\mathrm{d}x\right)^2}
& =\frac{\frac{(N+3\alpha+1)
(N+2\alpha)
\Gamma\left(\frac{N+2\alpha}
{\alpha+1}\right)}
{4(2b)^{\frac{N+2\alpha}{\alpha+1}}}
\frac{\Gamma
\left(\frac{N+2\alpha}
{\alpha+1}+2\right)}
{(2b)^{\frac{N+2\alpha}{\alpha+1}+2}}}
{\left[\frac{\Gamma
\left(\frac{N+2\alpha}
{\alpha+1}+1\right)}
{(2b)^{\frac{N+2\alpha}{\alpha+1}+1}}
\right]^2}
\\& =\frac{(N+3\alpha+1)
(N+2\alpha)
\left(\frac{N+2\alpha}
{\alpha+1}+1\right)
\frac{N+2\alpha}{\alpha+1}}
{4\left(\frac{N+2\alpha}
{\alpha+1}\right)^2}
\\&=\frac{(N+3\alpha+1)^2}{4}.
\end{align*}
Now, the proof of Theorem \ref{thmwsupsc-2}
is thereby completed.
\end{proof}

\subsection{The case $N=4$ and $-1<\alpha\le0$: proof of Theorem \ref{thmwsupscN4}}
\label{sect-3.3}

\noindent In this subsection, we prove Theorem \ref{thmwsupscN4}.

Let $\mathbb{H}$ be a Hilbert functional space, and denote
\begin{align}
\mu&=\inf_{u\in\mathbb{H}
\setminus\{0\}}
\frac{\sqrt{H(u)U(u)}}{P(u)},
\label{11.12-1}\\
\mu_+&=\inf_{u\in\mathbb{H}
\setminus\{0\}}
\frac{H(u)+U(u)}{P(u)},
\label{11.12-2}
\end{align}
where $H,U$ and $P$ are continuous and positive definite quadratic forms on $\mathbb{H}\setminus\{0\}$, and $\mu,\mu_+>0$ are constants. Without loss of generality, it is assumed that $\mathbb{H}$ is scaling-invariant for variables, and $H,U$ and $P$ are homogeneous with respect to the scaling. Let $u_\lambda(x)=u(\lambda x)\in\mathbb{H}\setminus\{0\}$ for all $\lambda>0$, there exists $(\alpha_1,\alpha_2,\alpha_3)
\in\mathbb{R}^3$ such that
\begin{equation}\label{si}
H(u_\lambda)=\lambda^{\alpha_1}H(u),\ \
U(u_\lambda)=\lambda^{\alpha_2}U(u),\ \
P(u_\lambda)=\lambda^{\alpha_3}P(u).
\end{equation}
In order to obtain \eqref{11.12-1}, there must hold  $\alpha_1+\alpha_2=2\alpha_3$. Inspired by Huang and Ye in \cite{Huang25}, the key point in the proof of Theorem \ref{thmwsupscN4} is to establish the equivalence between \eqref{11.12-1} and \eqref{11.12-2}, and we provide the proof of it in the following.

\begin{prop}\label{prop-1}
Under invariance of scaling, the inequality \eqref{11.12-1} is equivalent to the inequality \eqref{11.12-2}.
\end{prop}

\begin{proof}[\rm\textbf{Proof}]
We divide the proof into three steps.

\emph{First, we obtain that $\mu_+=2\mu$, where $\mu_+,\mu$ are defined by \eqref{11.12-1} and \eqref{11.12-2} respectively.} One hand, using Cauchy-Schwarz inequality, we see that
\begin{equation}\label{11.13-1}
\mu_+=\inf_{u\in\mathbb{H}\setminus\{0\}}
\frac{H(u)+U(u)}{P(u)}
\ge \inf_{u\in\mathbb{H}\setminus\{0\}}
\frac{2\sqrt{H(u)U(u)}}{P(u)}=2\mu.
\end{equation}
On the other hand, for all $u\in\mathbb{H}\setminus\{0\}$ and $\lambda>0$, there exists $u_\lambda\in\mathbb{H}\setminus\{0\}$ such that
\begin{equation*}
\mu_+=\inf_{u\in\mathbb{H}\setminus\{0\}}
\frac{H(u)+U(u)}{P(u)}
\le\frac{H(u_\lambda)+U(u_\lambda)}
{P(u_\lambda)}
=\frac{\lambda^{\alpha_1-\alpha_3}H(u)
+\lambda^{\alpha_2-\alpha_3}U(u)}
{P(u)},
\end{equation*}
this together with $\alpha_1+\alpha_2=2\alpha_3$ deduces that, for all $\lambda>0$,
\begin{equation*}
\mu_+\le
\frac{\lambda^{\alpha_1-\alpha_3}H(u)
+\lambda^{-(\alpha_1-\alpha_3)}U(u)}
{P(u)},
\end{equation*}
which derives that
\begin{equation*}
\mu_+\le\inf_{\lambda>0}
\frac{\lambda^{\alpha_1-\alpha_3}H(u)
+\lambda^{-(\alpha_1-\alpha_3)}U(u)}
{P(u)}
=\frac{2\sqrt{H(u)U(u)}}{P(u)}.
\end{equation*}
Due to the arbitrariness of $u$,
\begin{equation}\label{11.13-2}
\mu_+\le2
\inf_{u\in\mathbb{H}\setminus\{0\}}
\frac{\sqrt{H(u)U(u)}}{P(u)}
=2\mu.
\end{equation}
Combining \eqref{11.13-1} and \eqref{11.13-2}, we see that
\begin{equation*}
\mu_+=2\mu.
\end{equation*}
Therefore, \eqref{11.12-1} and \eqref{11.12-2} become
\begin{align*}
\mu&=\inf_{u\in\mathbb{H}\setminus\{0\}}
\frac{\sqrt{H(u)U(u)}}{P(u)},\\
2\mu&=\inf_{u\in\mathbb{H}\setminus\{0\}}
\frac{H(u)+U(u)}{P(u)}.
\end{align*}
This implies that
\begin{align}
H(u)U(u)&\ge \mu^2P(u)^2,\label{11.13-4}\\
H(u)+U(u)&\geq 2\mu P(u),\label{11.10-1}
\end{align}
and $\mu^2$ and $2\mu$ are the sharp constant of \eqref{11.13-4} and \eqref{11.10-1} respectively. In view of this, let
\begin{align*}
\Sigma_1
&:=\left\{u\in\mathbb{H}\setminus\{0\}:\mu=
\frac{\sqrt{H(u)U(u)}}{P(u)}\right\},\\
\Sigma_2
&:=\left\{u\in\mathbb{H}\setminus\{0\}:2\mu=
\frac{H(u)+U(u)}{P(u)}\right\}
\end{align*}
be the set of extremal functions of \eqref{11.13-4} and \eqref{11.10-1} respectively.

\emph{Second, we verify that the extremal function of the inequality \eqref{11.10-1} is that of the inequality \eqref{11.13-4}.} Let $u_\star\in\Sigma_2$ be the extremal function of \eqref{11.10-1}, applying Cauchy-Schwarz inequality and \eqref{11.13-4},
\begin{equation*}
2\mu=\frac{H(u_\star)+U(u_\star)}
{P(u_\star)}
\ge\frac{2\sqrt{H(u_\star)U(u_\star)}}
{P(u_\star)}
\ge2\mu,
\end{equation*}
that is,
\begin{equation*}
\mu=\frac{\sqrt{H(u_\star)U(u_\star)}}
{P(u_\star)},
\end{equation*}
which means that $u_\star\in\Sigma_1$ is the extremal function of \eqref{11.13-4}. Furthermore, from \eqref{si}, it is obvious that $u_{\star,\lambda}(x)=u_\star(\lambda x)\in\Sigma_1$ for all $\lambda>0$.

\emph{Last, we show that the extremal function of the inequality \eqref{11.13-4} is that of the inequality \eqref{11.10-1} up to a scaling.} Mathematically, for any $u\in\Sigma_1$, there exists $\lambda>0$ such that  $u_\lambda\in\Sigma_2$. Let $u_\diamond\in\Sigma_1$ be the extremal function of \eqref{11.13-4}, and $u_{\diamond,\lambda}
=u_\diamond(\lambda x)$ for $\lambda^{\alpha_1-\alpha_3}
=\sqrt{\frac{U(u_\diamond)}
{H(u_\diamond)}}$, by \eqref{11.10-1},
\begin{align*}
2\mu
\le\frac{H(u_{\diamond,\lambda})
+U(u_{\diamond,\lambda})}
{P(u_{\diamond,\lambda})}
=\frac{\lambda^{\alpha_1-\alpha_3}
H(u_\diamond)
+\lambda^{-(\alpha_1-\alpha_3)}
U(u_\diamond)}
{P(u_\diamond)}
=\frac{2\sqrt{H(u_\diamond)
U(u_\diamond)}}
{P(u_\diamond)}
=2\mu,
\end{align*}
that is,
\begin{equation*}
2\mu
=\frac{H(u_{\diamond,\lambda})
+U(u_{\diamond,\lambda})}
{P(u_{\diamond,\lambda})},
\end{equation*}
which implies that $u_{\diamond,\lambda}\in\Sigma_2$ is the extremal function of \eqref{11.10-1}.

In summary, under invariance of scaling, the inequality \eqref{11.13-4} is equivalent to the inequality \eqref{11.10-1}.
\end{proof}

Therefore, from Proposition \ref{prop-1}, to obtain the sharp form of the inequality \eqref{wsUPN4}, it suffices to consider the following inequality
\begin{equation*}
\int_{\mathbb{R}^4}
\frac{|\Delta u|^2}
{|x|^{2\alpha}}
\mathrm{d}x
+C^2\int_{\mathbb{R}^4}\left|\nabla u\right|^2\mathrm{d}x
\ge C(3\alpha+5)
\int_{\mathbb{R}^4}
\frac{\left|\nabla u\right|^2}
{|x|^{\alpha+1}}
\mathrm{d}x,
\end{equation*}
when choosing
\begin{equation*}
H(u)=\int_{\mathbb{R}^4}
\frac{|\Delta u|^2}
{|x|^{2\alpha}}
\mathrm{d}x,
\ \
U(u)=C^2\int_{\mathbb{R}^4}
\left|\nabla u\right|^2\mathrm{d}x,
\ \
P(u)=C\int_{\mathbb{R}^4}
\frac{\left|\nabla u\right|^2}
{|x|^{\alpha+1}}
\mathrm{d}x,
\ \
\mu=\frac{3\alpha+5}{2}
\end{equation*}
in \eqref{11.13-4} and \eqref{11.10-1}, for any $C>0$. Inspired by the work of Deng and Tian \cite{Deng23}, applying the change of variable $r\mapsto r^{\alpha+1}$, we indeed consider the following form
\begin{equation}\label{wsUPN4b}
\int_{\mathbb{R}^4}
\frac{|\Delta u|^2}
{|x|^{2\alpha}}
\mathrm{d}x
+q^{-2}\int_{\mathbb{R}^4}
\left|\nabla u\right|^2\mathrm{d}x
\ge q^{-1}(3\alpha+5)
\int_{\mathbb{R}^4}
\frac{\left|\nabla u\right|^2}
{|x|^{\alpha+1}}
\mathrm{d}x,
\end{equation}
where $q=\frac{1}{\alpha+1}$. Then, in order to prove Theorem \ref{thmwsupscN4}, based on the argument of Proposition \ref{prop-1}, it suffices to verify that the inequality \eqref{wsUPN4b} holds.

Let $u\in C^\infty_c(\mathbb{R}^4)$, and
\begin{equation*}
\mathcal{J}_{\alpha}(u):
=\int_{\mathbb{R}^4}
\frac{|\Delta u|^2}{|x|^{2\alpha}}\mathrm{d}x
+q^{-2}\int_{\mathbb{R}^4}\left|\nabla u\right|^2\mathrm{d}x
-q^{-1}(3\alpha+5)\int_{\mathbb{R}^4}
\frac{\left|\nabla u\right|^2}{|x|^{\alpha+1}}
\mathrm{d}x,
\end{equation*}
where $q:=\frac{1}{\alpha+1}$.
Using the spherical harmonic decomposition
and Lemma \ref{lemhduukd}, we have
\begin{equation*}
\mathcal{J}_{\alpha}(u)
=\sum\limits_{k=0}^\infty \mathcal{J}_{\alpha,k}(v_k),
\end{equation*}
where
\begin{align*}
\mathcal{J}_{\alpha,k}(v)
&=\int_0^\infty |v''|^2r^{2k-2\alpha+3}
\mathrm{d}r
+(2\alpha+1)(2k+3)\int_0^\infty |v'|^2
r^{2k-2\alpha+1}
\mathrm{d}r
\nonumber\\&\quad
+q^{-2}\int_0^\infty |v'|^2r^{2k+3}\mathrm{d}r
-q^{-1}(3\alpha+5)\int_0^\infty |v'|^2r^{2k-\alpha+2}\mathrm{d}r
\nonumber\\&\quad
-q^{-1}(3\alpha+5)(\alpha+1)k
\int_0^\infty |v|^2r^{2k-\alpha}
\mathrm{d}r.
\end{align*}

For each $v\in C_{c}^{2}([0,\infty))$, let us make the change
\begin{equation*}
v(r)=w(s),
\ \
\mbox{where}
\
r=s^q,
\
q=\frac{1}{\alpha+1}.
\end{equation*}
A direct calculation shows that
\begin{align*}
v'(r)& =\frac{\mathrm{d}v(r)}{\mathrm{d}r}
=\frac{\mathrm{d}w(s)}{\mathrm{d}s}
\frac{\mathrm{d}s}{\mathrm{d}r}
=w'(s)q^{-1}r^{q^{-1}-1}
=q^{-1}w'(s)s^{1-q},
\nonumber\\
v''(r)& =\frac{\mathrm{d}v'(r)}{\mathrm{d}r}
=\frac{\mathrm{d}
\left(q^{-1}w'(s)s^{1-q}\right)}
{\mathrm{d}s}
\frac{\mathrm{d}s}{\mathrm{d}r}
=q^{-2}w''(s)s^{2-2q}
+q^{-2}(1-q)w'(s)s^{1-2q},
\end{align*}
then using integration by parts, we have
\begin{align*}
\int_0^\infty |v''|^2r^{2k-2\alpha+3}\mathrm{d}r
&=q^{-3}\int_0^\infty |w''|^2
s^{q(2k-2\alpha)+3}\mathrm{d}s
\\&\quad -q^{-3}(1-q)[q(2k-2\alpha+1)+1]
\int_0^\infty |w'|^2
s^{q(2k-2\alpha)+1}\mathrm{d}s,
\end{align*}
and
\begin{align*}
\int_0^\infty |v'|^2r^{2k-2\alpha+1}\mathrm{d}r
&=q^{-1}\int_0^\infty |w'|^2
s^{q(2k-2\alpha)+1}\mathrm{d}s,
\\
\int_0^\infty |v'|^2r^{2k+3}\mathrm{d}r
&=q^{-1}\int_0^\infty |w'|^2
s^{q(2k+2)+1}\mathrm{d}s,
\\
\int_0^\infty |v'|^2r^{2k-\alpha+2}\mathrm{d}r
&=q^{-1}\int_0^\infty |w'|^2
s^{q(2k-\alpha+1)+1}\mathrm{d}s,
\\
\int_0^\infty |v|^2r^{2k-\alpha}\mathrm{d}r
&=q\int_0^\infty |w|^2
s^{q(2k-\alpha+1)-1}\mathrm{d}s.
\end{align*}
Therefore,
\begin{align}\label{10.22-1b}
\mathcal{J}_{\alpha,k}(v)
&=q^{-3}\int_0^\infty |w''|^2
s^{q(2k-2\alpha)+3}\mathrm{d}s
-q^{-3}(1-q)[q(2k-2\alpha+1)+1]
\int_0^\infty |w'|^2
s^{q(2k-2\alpha)+1}\mathrm{d}s
\nonumber\\&\quad
+q^{-1}(2\alpha+1)(2k+3)\int_0^\infty |w'|^2
s^{q(2k-2\alpha)+1}\mathrm{d}s
+q^{-3}\int_0^\infty |w'|^2
s^{q(2k+2)+1}\mathrm{d}s
\nonumber\\&\quad
-q^{-2}(3\alpha+5)\int_0^\infty |w'|^2
s^{q(2k-\alpha+1)+1}\mathrm{d}s
-(3\alpha+5)(\alpha+1)k\int_0^\infty |w|^2
s^{q(2k-\alpha+1)-1}\mathrm{d}s
\nonumber\\&=
q^{-3}\int_0^\infty |w''+w'|^2
s^{q(2k-2\alpha)+3}\mathrm{d}s
+q^{-3}[q(2k-2\alpha)+3]
\int_0^\infty |w'|^2
s^{q(2k-2\alpha)+2}\mathrm{d}s
\nonumber\\&\quad
+\left\{q^{-1}(2\alpha+1)(2k+3)
-q^{-3}(1-q)[q(2k-2\alpha+1)+1]\right\}
\int_0^\infty |w'|^2
s^{q(2k-2\alpha)+1}\mathrm{d}s
\nonumber\\&\quad
-q^{-2}(3\alpha+5)\int_0^\infty |w'|^2
s^{q(2k-\alpha+1)+1}\mathrm{d}s -(3\alpha+5)(\alpha+1)k\int_0^\infty |w|^2
s^{q(2k-\alpha+1)-1}\mathrm{d}s
\nonumber\\&=q^{-3}\int_0^\infty |w''+w'|^2
s^{q(2k-2\alpha)+3}\mathrm{d}s
+2q^{-2}(k-\alpha-1)
\int_0^\infty |w'|^2
s^{q(2k-2\alpha)+2}\mathrm{d}s
\nonumber\\&\quad
+\left\{q^{-1}(2\alpha+1)(2k+3)
-q^{-3}(1-q)[q(2k-2\alpha+1)+1]\right\}
\int_0^\infty |w'|^2
s^{q(2k-2\alpha)+1}\mathrm{d}s
\nonumber\\&\quad
-(3\alpha+5)(\alpha+1)k\int_0^\infty |w|^2
s^{q(2k-2\alpha)}\mathrm{d}s
\nonumber\\:\!&=\mathcal{K}_{\alpha,k}(w),
\end{align}
where the second equality and the third equality follow from
\begin{align*}
&q^{-3}\int_0^\infty |w''|^2
s^{q(2k-2\alpha)+3}\mathrm{d}s
+q^{-3}\int_0^\infty |w'|^2
s^{q(2k+2)+1}\mathrm{d}s
\\&\quad
=
q^{-3}\int_0^\infty |w''|^2
s^{q(2k-2\alpha)+3}\mathrm{d}s
+q^{-3}\int_0^\infty |w'|^2
s^{q(2k-2\alpha)+3}\mathrm{d}s
\\&\quad
=q^{-3}\int_0^\infty |w''+w'|^2
s^{q(2k-2\alpha)+3}\mathrm{d}s
-2q^{-3}\int_0^\infty w''w'
s^{q(2k-2\alpha)+3}\mathrm{d}s
\\&\quad
=q^{-3}\int_0^\infty |w''+w'|^2
s^{q(2k-2\alpha)+3}\mathrm{d}s
+q^{-3}[q(2k-2\alpha)+3]\int_0^\infty |w'|^2
s^{q(2k-2\alpha)+2}\mathrm{d}s,
\end{align*}
and
\begin{align*}
&q^{-3}[q(2k-2\alpha)+3]
\int_0^\infty |w'|^2
s^{q(2k-2\alpha)+2}\mathrm{d}s
-q^{-2}(3\alpha+5)\int_0^\infty |w'|^2
s^{q(2k-\alpha+1)+1}\mathrm{d}s
\\&\quad
=
q^{-3}[q(2k-2\alpha)+3]
\int_0^\infty |w'|^2
s^{q(2k-2\alpha)+2}\mathrm{d}s
-q^{-2}(3\alpha+5)
\int_0^\infty |w'|^2
s^{q(2k-2\alpha)+2}\mathrm{d}s
\\&\quad=2q^{-2}(k-\alpha-1)
\int_0^\infty |w'|^2
s^{q(2k-2\alpha)+2}\mathrm{d}s,
\end{align*}
respectively, with the aid of $q=\frac{1}{\alpha+1}$.

In view of this, to show $\mathcal{J}_{\alpha,k}(v)\ge0$ is equivalent to show $\mathcal{K}_{\alpha,k}(w)\ge0$.

\begin{lem}\label{k0k2}
Assume that $-1<\alpha\le0$. Let $v\in C_{c}^{2}([0,\infty))$,
then $\mathcal{K}_{\alpha,k}(w)\geq 0$ for all $k\geq 0$.
\end{lem}

\begin{proof}[\rm\textbf{Proof}]
\emph{We first consider $\mathcal{K}_{\alpha,0}(w)$.} Let $w(s)=e^{-s}\xi(s)$, there hold $w'=e^{-s}(\xi'-\xi)$ and $w''+w'=e^{-s}(\xi''-\xi')$. Let us denote $\varrho(s)=\xi'(s)-\xi(s)$, then $w'=e^{-s}\varrho$ and $w''+w'=e^{-s}\varrho'$, and thus
\begin{align}\label{k=0}
\mathcal{K}_{\alpha,0}(w)
&=q^{-3}\int_0^\infty |w''+w'|^2
s^{-2q\alpha+3}\mathrm{d}s
-2q^{-2}(\alpha+1)
\int_0^\infty |w'|^2
s^{-2q\alpha+2}\mathrm{d}s
\nonumber\\&\quad +\left\{3 q^{-1}(2\alpha+1)
-q^{-3}(1-q)[q(1-2\alpha)+1]\right\}
\int_0^\infty |w'|^2
s^{-2q\alpha+1}\mathrm{d}s
\nonumber\\
&=q^{-3}\int_0^\infty |\varrho'|^2
e^{-2s}s^{-2q\alpha+3}\mathrm{d}s
-2q^{-2}(\alpha+1)
\int_0^\infty |\varrho|^2e^{-2s}
s^{-2q\alpha+2}\mathrm{d}s
\nonumber\\&\quad +\left\{3 q^{-1}(2\alpha+1)
-q^{-3}(1-q)[q(1-2\alpha)+1]\right\}
\int_0^\infty |\varrho|^2e^{-2s}
s^{-2q\alpha+1}\mathrm{d}s
\nonumber\\
&=q^{-3}\int_0^\infty (s\varrho'-\varrho)^2
e^{-2s}s^{-2q\alpha+1}\mathrm{d}s
\geq 0.
\end{align}

\emph{Next, we consider $\mathcal{K}_{\alpha,k}(w)$ for $k\geq 2$.} Taking $\theta=q(2k-2\alpha)+2$ in the Hardy inequality \eqref{HI},
\begin{align}\label{HIi}
&2q^{-2}(k-\alpha-1)
\int_0^\infty |w'|^2
s^{q(2k-2\alpha)+2}\mathrm{d}s
\nonumber\\&\quad\geq \frac{q^{-2}(k-\alpha-1)[q(2k-2\alpha)+1]^2}{2}
\int_0^\infty |w|^2
s^{q(2k-2\alpha)}\mathrm{d}s
\nonumber\\
&\quad\geq (3\alpha+5)(\alpha+1)k\int_0^\infty |w|^2
s^{q(2k-2\alpha)}\mathrm{d}s.
\end{align}
where the second inequality follows from
\begin{align*}
\frac{q^{-2}(k-\alpha-1)
[q(2k-2\alpha)+1]^2}{2}
&=\frac{(k-\alpha-1)(2k-\alpha+1)^2}{2}
\\&\geq
\frac{k(2k+1)^2}{4}
\\&\geq
\frac{25k}{4}\geq (3\alpha+5)(\alpha+1)k,
\end{align*}
for all $-1<\alpha\le0$. Furthermore,
\begin{align}\label{k2}
& q^{-1}(2\alpha+1)(2k+3)-q^{-3}(1-q)
[q(2k-2\alpha+1)+1]
\nonumber\\
&\quad=2k(\alpha+1)^2
+(\alpha+1)^2(\alpha+3)\ge0.
\end{align}
Inserting \eqref{HIi} and \eqref{k2} into \eqref{10.22-1b}, we have $\mathcal{K}_{\alpha,k}(w)\geq 0$ for $k\geq 2$.

\emph{Lastly, it remains to consider $\mathcal{K}_{\alpha,1}(w)$.} Let us recall from \eqref{10.22-1b} that
\begin{align*}
\mathcal{K}_{\alpha,1}(w)
&=q^{-3}\int_0^\infty |w''+w'|^2
s^{q(2-2\alpha)+3}\mathrm{d}s
-2q^{-2}\alpha
\int_0^\infty |w'|^2
s^{q(2-2\alpha)+2}\mathrm{d}s
\nonumber\\&\quad
+\left\{5q^{-1}(2\alpha+1)
-q^{-3}(1-q)[q(2-2\alpha+1)+1]\right\}
\int_0^\infty |w'|^2
s^{q(2-2\alpha)+1}\mathrm{d}s
\nonumber\\&\quad
-(3\alpha+5)(\alpha+1)\int_0^\infty |w|^2
s^{q(2-2\alpha)}\mathrm{d}s
\\&=q^{-3}\int_0^\infty |w''+w'|^2
s^{q(2-2\alpha)+3}\mathrm{d}s
-2q^{-2}\alpha
\int_0^\infty |w'|^2
s^{q(2-2\alpha)+2}\mathrm{d}s
\nonumber\\&\quad +(\alpha+1)^2(\alpha+5)
\int_0^\infty |w'|^2
s^{q(2-2\alpha)+1}\mathrm{d}s -(3\alpha+5)(\alpha+1)\int_0^\infty |w|^2
s^{q(2-2\alpha)}\mathrm{d}s,
\end{align*}
where the last equality follows from \eqref{k2} (with $k=1$). Then due to $-1<\alpha\leq 0$, choosing $\theta=\frac{4}{\alpha+1}$ in the Hardy inequality \eqref{HI}, we have
\begin{align*}
\mathcal{K}_{\alpha,1}(w)
&=(\alpha+1)^3\int_0^\infty |w''+w'|^2
s^{\frac{4}{\alpha+1}+1}\mathrm{d}s
-2\alpha(\alpha+1)^2
\int_0^\infty |w'|^2
s^{\frac{4}{\alpha+1}}\mathrm{d}s
\nonumber\\&\quad +(\alpha+1)^2(\alpha+5)
\int_0^\infty |w'|^2
s^{\frac{4}{\alpha+1}-1}\mathrm{d}s
-(3\alpha+5)(\alpha+1)\int_0^\infty |w|^2
s^{\frac{4}{\alpha+1}-2}\mathrm{d}s
\nonumber\\&\geq (\alpha+1)^3\int_0^\infty |w''+w'|^2
s^{\frac{4}{\alpha+1}+1}\mathrm{d}s +(\alpha+1)^2(\alpha+5)\int_0^\infty |w'|^2
s^{\frac{4}{\alpha+1}-1}\mathrm{d}s
\nonumber\\&\quad -\left[\frac{\alpha(3-\alpha)^2}{2}
+(3\alpha+5)(\alpha+1)\right]
\int_0^\infty |w|^2
s^{\frac{4}{\alpha+1}-2}\mathrm{d}s
\nonumber\\& \geq 4(\alpha+1)\int_0^\infty |w'+w|^2
s^{\frac{4}{\alpha+1}-1}\mathrm{d}s +(\alpha+1)^2(\alpha+5)\int_0^\infty |w'|^2
s^{\frac{4}{\alpha+1}-1}\mathrm{d}s
\nonumber\\&\quad
-\left[\frac{\alpha(3-\alpha)^2}{2}
+(3\alpha+5)(\alpha+1)\right]
\int_0^\infty |w|^2
s^{\frac{4}{\alpha+1}-2}\mathrm{d}s
\nonumber\\&= 4(\alpha+1)\int_0^\infty |\xi'|^2e^{-2s}
s^{\frac{4}{\alpha+1}-1}\mathrm{d}s
+(\alpha+1)^2(\alpha+5)\int_0^\infty |\xi'-\xi|^2e^{-2s}
s^{\frac{4}{\alpha+1}-1}\mathrm{d}s
\nonumber\\&\quad
-\left[\frac{\alpha(3-\alpha)^2}{2}
+(3\alpha+5)(\alpha+1)\right]
\int_0^\infty |\xi|^2e^{-2s}
s^{\frac{4}{\alpha+1}-2}\mathrm{d}s
\nonumber\\&= (\alpha+1)(\alpha+3)^2\int_0^\infty |\xi'|^2e^{-2s}
s^{\frac{4}{\alpha+1}-1}\mathrm{d}s
- (\alpha+1)^2(\alpha+5)\int_0^\infty |\xi|^2e^{-2s}
s^{\frac{4}{\alpha+1}-1}\mathrm{d}s
\nonumber\\&\quad
+\left[(\alpha+1)(\alpha+5)(3-\alpha)
-\frac{\alpha(3-\alpha)^2}{2}
-(3\alpha+5)(\alpha+1)\right]
\int_0^\infty |\xi|^2e^{-2s}
s^{\frac{4}{\alpha+1}-2}\mathrm{d}s
\nonumber\\:\!&=\mathcal{R}_{\alpha,1}(\xi),
\end{align*}
here the change of variable $w(s)=e^{-s}\xi(s)$ is required to obtain the second equality. Furthermore, applying the Hardy inequality \eqref{HIk1} with $\theta=\frac{4}{\alpha+1}-1$ and $\kappa=\frac{\alpha+1}{\alpha+3}$, we have
\begin{align}\label{10.22-1b3k1f}
\int_0^\infty |\xi'|^2e^{-2s}
s^{\frac{4}{\alpha+1}-1}\mathrm{d}s
& \geq
 -\frac{3-\alpha}{\alpha+3}\int_0^\infty |\xi|^2 e^{-2s} s^{\frac{4}{\alpha+1}-2} \mathrm{d}s
\nonumber\\&\quad +\frac{(\alpha+1)(\alpha+5)}{(\alpha+3)^2}
\int_0^\infty |\xi|^2 e^{-2s} s^{\frac{4}{\alpha+1}-1} \mathrm{d}s.
\end{align}
Inserting \eqref{10.22-1b3k1f} into the expression of $\mathcal{R}_{\alpha,1}(\xi)$, we obtain
\begin{align*}
\mathcal{R}_{\alpha,1}(\xi)&\geq  \left[(\alpha+1)(\alpha+5)(3-\alpha)
-\frac{\alpha(3-\alpha)^2}{2}
-(3\alpha+5)(\alpha+1)
-(\alpha+1)(\alpha+3)(3-\alpha)\right]
\nonumber\\&\quad\times\int_0^\infty |\xi|^2e^{-2s}
s^{\frac{4}{\alpha+1}-2}\mathrm{d}s
\nonumber\\&=\left[(\alpha+1)(1-5\alpha)
-\frac{\alpha(3-\alpha)^2}{2}\right]
\int_0^\infty |\xi|^2e^{-2s}
s^{\frac{4}{\alpha+1}-2}\mathrm{d}s,
\end{align*}
hence $\mathcal{K}_{\alpha,1}(w)\geq \mathcal{R}_{\alpha,1}(\xi)\geq 0$ due to $(\alpha+1)(1-5\alpha)
-\frac{\alpha(3-\alpha)^2}{2}\geq 0$ for all $-1<\alpha\leq 0$.
\end{proof}

\begin{lem}\label{k=1}
Assume that $-1<\alpha\leq 0$. For all $u\in \mathcal{H}^2_{\alpha,0}(\mathbb{R}^4)$, then $\mathcal{J}_{\alpha}(u)\geq0$. Furthermore, $\mathcal{J}_{\alpha}(u)=0$ if and only if $u(x)=A(1+|x|^{\alpha+1})
\exp(-|x|^{\alpha+1})$ for all $A\in\mathbb{R}$.
\end{lem}

\begin{proof}[\rm\textbf{Proof}]
From Lemma \ref{k0k2}, we see that  $\mathcal{J}_{\alpha}(u)\geq 0$ for all $u\in C^2_c(\mathbb{R}^4)$, and $\mathcal{J}_{\alpha}(u)\geq 0$ also holds for all $u\in \mathcal{H}^2_{\alpha,0}(\mathbb{R}^4)$ by the density argument.

Moreover, by the above analysis and the non-attainability of the Hardy inequality \eqref{HI} for any nontrivial function, we see that to reach $\mathcal{J}_{\alpha}(u)
=\sum_{k=0}^\infty \mathcal{J}_{\alpha,k}(v_k)=0$, the only choice is $v_{k}\equiv 0$ for all $k\geq1$; and the equality holds in \eqref{k=0}. This means that $s\varrho^{\prime}\equiv\varrho$, that is $\xi^{\prime}-\xi=\varrho=As$, which indicates $\xi(s)=A(1+s)$, then $u(x)
=v_0(r)=A(1+r^{\alpha+1})
\exp(-r^{\alpha+1})$ with $A\in\mathbb{R}$.
\end{proof}

We finish this section by proving Theorem \ref{thmwsupscN4}.

\begin{proof}
[\rm\textbf{Proof of Theorem \ref{thmwsupscN4}}]
By Proposition \ref{prop-1}, the HyUP inequality \eqref{wsUPN4} is equivalent to \eqref{wsUPN4b} up to a scaling. From Lemma \ref{k=1}, we know that inequality \eqref{wsUPN4b} is sharp and the equality holds if and only if $u(x)=A(1+|x|^{\alpha+1})
\exp(-|x|^{\alpha+1})$ with $A\in\mathbb{R}$. Thus, we deduce that the HyUP inequality \eqref{wsUPN4} holds for any $u\in \mathcal{H}^2_{\alpha,0}(\mathbb{R}^4)$,
where the constant $\frac{(3\alpha+5)^2}{4}$
is sharp and is attained if and only if
$u(x)=a(1+b|x|^{\alpha+1})
\exp(-b|x|^{\alpha+1})$ with $a\in\mathbb{R}$ and $b>0$.
\end{proof}

\section*{Declarations}

\subsection*{Funding}
\noindent This paper was supported by the National Natural Science Foundation of China (No. 12371120).

\subsection*{Data availability statement}
\noindent No data was used for the research described in the article.

\subsection*{Conflict of interest}
\noindent The authors declare no conflict of interest.

\end{document}